\documentclass{article}
\usepackage{latexsym}
\usepackage{graphics}
\usepackage{amsmath}
\usepackage{amsfonts}
\begin{document}
\bibliographystyle{plain}
\title{On the decomposition of Global Conformal Invariants II}
\author{Spyros Alexakis}
\maketitle
\newtheorem{proposition}{Proposition}
\newtheorem{theorem}{Theorem}
\newtheorem{lemma}{Lemma}
\newtheorem{conjecture}{Conjecture}
\newtheorem{observation}{Observation}
\newtheorem{formulation}{Formulation}
\newtheorem{definition}{Definition}
\newtheorem{corollary}{Corollary}

\begin{abstract}

\par This paper is a continuation of \cite{a:dgciI}, where we
 complete our partial proof of the Deser-Schwimmer conjecture
 on the structure of ``global conformal invariants''. Our
 theorem
deals with such invariants $P(g^n)$ that locally depend
only on the curvature tensor $R_{ijkl}$ (without covariant
 derivatives).

\par  In \cite{a:dgciI} we developed
 a powerful tool, the ``super divergence formula'' which
 applies to any Riemannian operator that always integrates to
 zero on compact manifolds. In particular, it applies to the
 operator $I_{g^n}(\phi)$ that measures the ``non-conformally
 invariant part'' of $P(g^n)$. This paper resolves the problem
 of using this information we have obtained on the structure
 of $I_{g^n}(\phi)$ to understand the structure of $P(g^n)$.
\end{abstract}

\section{Introduction}

\par We briefly recall the open problem that this paper and
 \cite{a:dgciI} address and the theorem that we will be
 completing here. Our objects of study are
{\it scalar Riemannian invariants} $P(g^n)$ of a Riemannian
 manifold $(M^n,g^n)$. These are polynomials in the components of
 the tensors $R_{ijkl},\dots ,\nabla^m_{r_1\dots r_m}R_{ijkl},
\dots$ and $g^{ij}$ (or, even more generally, in the variables
$\partial^k_{t_1\dots t_k}g_{ij}$, $det(g)^{-1}$), that are
 independent of the
coordinate system in which they are expressed, and also have a
weight $W$, meaning that under a re-scaling $g^n\rightarrow
t^2g^n$ they transform by $P(t^2g^n)=t^WP(g^n)$, $t\in
\mathbb{R}_{+}$. It is a classical result that such invariants are
linear combinations

\begin{equation}
\label{lincob}
P(g^n)=\Sigma_{l\in L} a_l C^l(g^n)
\end{equation}

of complete contractions in the form:

\begin{equation}
\label{suntmhsh}
contr(\nabla^{m_1}_{r_1\dots r_{m_1}}R_{ijkl}\otimes
\dots\otimes\nabla^{m_s}_{t_1\dots t_{m_s}}R_{i'j'k'l'})
\end{equation}
each with weight $W$. We fix an {\it even} dimension $n$ once and
 for all, and we restrict attention to local scalar
invariants of weight $-n$. Due to the transformation of the
 volume form $dV_{e^{2\phi(x)}g^n}=e^{n\phi(x)}dV_{g^n}$
under general conformal re-scalings
$\hat{g}^n\rightarrow e^{2\phi(x)}g^n$, it follows that if
$P(g^n)$ has weight $-n$ then the quantity $\int_{M^n}
P(g^n)dV_{g^n}$ is scale-invariant for any compact orientable
Riemannian $(M^n,g^n)$.

\par The problem we are addressing is to find all Riemannian
scalar invariants of weight $-n$ for which the integral
$\int_{M^n}P(g^n)dV_{g^n}$ is invariant under conformal
re-scalings $\hat{g}^n=e^{2\phi(x)}g^n$ for
 any compact manifold
$(M^n,g^n)$ and any $\phi\in C^\infty (M^n)$. In other words, we
are assuming that for any $(M^n,g^n)$ and $\phi\in C^\infty (M^n)$
we must have:

\begin{equation}
\label{elie}
\int_{M^n}P(g^n)dV_{g^n}=\int_{M^n}P(\hat{g}^n)dV_{\hat{g}^n}
\end{equation}

Deser and Schwimmer, two physicists, conjectured the following in \cite{ds:gccaad}:

\begin{conjecture}[Deser-Schwimmer]
\label{conj} Suppose we have a Riemannian scalar $S(g^n)$ of
weight $-n$ for some even $n$. Suppose that for any compact
manifold $(M^n,g^n)$ the quantity

\begin{equation}
\label{donne} {\int}_{M^n}S(g^n)dV_{g^n}
\end{equation}
is invariant under any conformal change of metric
$\hat{g}^n(x)=
 e^{2\phi (x)} g^n(x)$. Then $P(g^n)$ must be a linear combination of
 three``obvious candidates'', namely:

\begin{equation}
\label{post} S(g^n)=W(g^n)+div_iT_i(g^n)+c\cdot
\operatorname{Pfaff}(R_{ijkl})
\end{equation}

\begin{enumerate}

\item{$W(g^n)$ is a scalar conformal invariant of weight
$-n$, ie
it satisfies $W(e^{2\phi(x)}g^n)\\ =e^{-n\phi(x)}W(g^n)$
for every $\phi\in C^\infty(M^n)$ and every $x\in M^n$.}

\item{$T_i(g^n)$ is a Riemannian vector field of weight
$-n+1$.
(Since for any compact $M^n$ we have $\int_{M^n}div_i
T_i(g^n)dV_{g^n}=0$.)}

\item{$\operatorname{Pfaff}(R_{ijkl})$ stands for the Pfaffian of
the curvature $R_{ijkl}$. (Since for any compact Riemannian
$(M^n,g^n)$ $\int_{M^n} \operatorname{Pfaff}(R_{ijkl})dV_{g^n}=
\frac{2^n\pi^{\frac{n}{2}}(\frac{n}{2}-1)!}{2(n-1)!}\chi(M^n)$.)}
\end{enumerate}
\end{conjecture}

\par In this paper we complete our partial confirmation of this conjecture.
 We restrict our attention to Riemannian
 scalars $P(g^n)$  that are linear combinations
\begin{equation}
\label{bombaster} \Sigma_{l\in L} a_l C^l(g^n)
\end{equation}
 of complete
contractions of weight $-n$, each $C^l(g^n)$ in the form:

\begin{equation}
\label{noderivatives}
contr(R_{i_1j_1k_1l_1}\otimes\dots\otimes
R_{i_{\frac{n}{2}}j_{\frac{n}{2}}k_{\frac{n}{2}}
l_{\frac{n}{2}}})
\end{equation}
(since we are not allowing derivatives on the factors
$R_{ijkl}$, the weight restriction forces each complete
contraction to have $\frac{n}{2}$ factors). The main theorem
that we show in \cite{a:dgciI} and in the present paper is:

\begin{theorem}
\label{weakt2}
Let us suppose that $P(g^n)$ is in the form (\ref{bombaster}),
where each $C^l(g^n)$ is in the form (\ref{noderivatives}),
 with $r=\frac{n}{2}$ factors. We also assume that
(\ref{elie}) holds for any Riemannian $(M^n,g^n)$ and
$\phi\in C^\infty (M^n)$.

Then, there exists a a scalar conformal invariant $W(g^n)$ of
 weight $-n$ that locally depends only on the Weyl tensor,
and also a constant $c$ so that:

\begin{equation}
\label{post2} S(g^n)=W(g^n)+c\cdot \operatorname{Pfaff}
(R_{ijkl})
\end{equation}
where $\operatorname{Pfaff}(R_{ijkl})$ stands for the Pfaffian of the
curvature $R_{ijkl}$.
\end{theorem}

\par We will recall two related results that were proven by entirely
different methods. In \cite{g:lierm} Gilkey
 considered the problem of finding all scalar invariants
$P(g^n)$ of weight
$-n$ for which $\int_{M^n}P(g^n)dV_{g^n}$ is constant for a given
compact orientable $M^n$ and {\it any} Riemannian
 metric $g^n$ over $M^n$. He then showed that:

\begin{theorem}[Gilkey]
\label{gilk}
Under the above assumptions, we have that $P(g^n)$ can be
written as:

\begin{equation}
\label{eqngilk}
P(g^n)=div_iT_i(g^n)+c\cdot \operatorname{Pfaff}(R_{ijkl})
\end{equation}
where $T_i(g^n)$ is an intrinsic vector field of weight
$-n+1$ and $\operatorname{Pfaff}(R_{ijkl})$ stands for the
 Pfaffian of the curvature tensor.
\end{theorem}
(see also \cite{m:pt} for an earlier form of this result).
Extending the methods in \cite{g:lierm}, Branson, Gilkey and
 Pohjanpelto showed in \cite{bgp:ilcfm} that:

\begin{theorem}[Branson-Gilkey-Pohjanpelto]
\label{BrGP}
\par Consider any local Riemannian invariant $P(g^n)$ of weight $-n$,
with the property that for any manifold $M^n$ and any locally conformally flat metric $h^n$,
$\int_{M^n}P(h^n)dV_{h^n}$  is invariant under conformal
re-scalings $\hat{h}^n=e^{2\phi(x)}h^n$ of the metric $h^n$.
 It then follows that in the locally conformally flat metric
 $h^n$ (for which the Weyl tensor vanishes), we can write out:

\begin{equation}
\label{mejico}
P(h^n)= div_i T_i(h^n)+c\cdot \operatorname{Pfaff}(R_{ijkl})
\end{equation}
where $T_i(h^n)$ is a vector field of weight $-n+1$ and
$\operatorname{Pfaff}(R_{ijkl})$ stands for the Pfaffian of
 the curvature tensor.
\end{theorem}

\par We have explained in \cite{a:dgciI} how resolving the whole of the
Deser-Schwimmer conjecture would have implications regarding the
structure of the so-called $Q$-curvature, and also for the study
of conformally compact Einstein manifolds, in particular regarding
the notions of the re-normalized volume and the conformal anomaly,
see also \cite{a:rcipem}, \cite{cqy:pc}, \cite{g:varccem},
\cite{gz:smcg}, \cite{gw:casoacc},\cite{hs:hwa}. Here, we briefly
recall the definition of $Q$-curvature.

\par $Q$-curvature is a Riemannian scalar invariant $Q^n(g^n)$
 constructed by Branson  for each even dimension $n$ (see \cite{b:fd}). In
 dimension 2 it is just the scalar curvature ($Q^2(g^2)=R$)
and in dimension 4 (where it has been extensively studied),
 it is in the form:

\begin{equation}
\label{tobecomple}
Q^4(g^4)= \frac{1}{12}\big{(} -\Delta R+\frac{1}{4}R^2-|E|^2\big{)}
\end{equation}
where $R$ is the scalar curvature and $E$ is the traceless
 Ricci  tensor.

\par In dimension $n$ $Q^n(g^n)$ has weight $-n$. Its two
 main properties are that $\int_{M^n}Q^n(g^n)dV_{g^n}$
 is invariant under conformal changes of $g^n$ and that
 under the re-scaling $g^n\rightarrow e^{2\phi(x)}g^n$, $Q^n(g^n)$
 enjoys the transformation law:

\begin{equation}
\label{komo}
Q^n(e^{2\phi(x)}g^n)(x)=e^{-n\phi(x)}[Q^n(g^n)+
P^{\frac{n}{2}}_{g^n}(\phi)](x)
\end{equation}
where $P^{\frac{n}{2}}_{g^n}(\phi)$ is a {\it conformally
 co-variant differential operator}, originally constructed in
 \cite{gjms:cipl}. Conformal co-variance means that its symbol
 has a nice transformation law under the conformal
re-scaling $\hat{g}^n=e^{2\phi(x)}g^n$, namely for every
$g^n$, $\phi,\psi\in C^\infty(M^n)$:

\begin{equation}
\label{covariance}
P^{\frac{n}{2}}_{e^{2\psi(x)}g^n}(\phi)=e^{-n\psi(x)}
P^{\frac{n}{2}}_{g^n}(\phi)
\end{equation}

\par The above transformation law has played an important role
in the analysis surrounding $Q$-curvature
(see \cite{cgy:ematcg}, \cite{b:gechofcg} for example). Moreover,
the particular form of $Q^4(g^4)$ and its relation to the
Chern-Gauss-Bonnet integrand has proven to be a valuable tool in
geometric and
 topological applications of $Q$-curvature in dimension 4,
see \cite{cqy:tccem}, \cite{q:rccem}. Therefore,  understanding
of the structure of $Q$-curvature in high
 dimensions would raise the question whether the powerful
 techniques employed in the study of $Q$-curvature in
dimension 4 can be extended to higher dimensions.

\section{Formulas and an outline of the proof.}

\par Throughout this paper we will be employing all the notational and
terminological conventions from \cite{a:dgciI}. We will also be
heavily using Theorem 2 in that paper and its two corollaries
regarding identities that hold ``formally'' or ``by
substitution'', see also \cite{beg:itccg}, \cite{e:ntrm},
\cite{w:cg}.
\newline

\par We recall that $P(g^n)$ satisfies (\ref{elie}).
 In \cite{a:dgciI} we defined an operator $I_{g^n}(\phi)$ as:

\begin{equation}
\label{nostos}
I_{g^n}(\phi)=e^{n\phi(x)}P(e^{2\phi(x)}g^n)-P(g^n)
\end{equation}
 which has weight $-n$ and the fundamental property that:

\begin{equation}
\label{elie2}
\int_{M^n}I_{g^n}(\phi)dV_{g^n}=0
\end{equation}
for every compact Riemannian $(M^n,g^n)$.
\newline

\par As our tool for this paper will be the super divergence
 formula for $I_{g^n}(\phi)$, it is necessary to
write out $P(g^n)$ in such a way so that we can ``recover''
the non-conformally invariant part of $P(g^n)$ from the
expression of $I_{g^n}(\phi)$.  As an illustration of the
 difficulty that we are forced to address, we
suppose that we write out $P(g^n)$ as a linear combination of
contractions in the form (\ref{noderivatives}). But then, given
the transformation law for the curvature tensor, it is not obvious
how to reconstruct $P(g^n)$ if we are given $I_{g^n}(\phi)$.
\newline

\par In order to overcome this difficulty, we recall
the Schouten tensor as a trace-adjustment of Ricci curvature:

\begin{equation}
\label{WSch}
P_{\alpha\beta}=\frac{1}{n-2}[{Ric_{\alpha\beta}}-\frac{R}{2(n-1)}
g^n_{\alpha\beta}]
\end{equation}

Where $Ric_{\alpha\beta}$ stands for Ricci curvature and $R$
stands for scalar curvature. We then have the well-known decomposition
of the curvature tensor:

\begin{equation}
\label{Weyl} R_{ijkl}=
W_{ijkl}+[P_{jk}g^n_{il}+P_{il}g^n_{jk}-P_{jl}g^n_{ik}-P_{ik}g^n_{jl}]
\end{equation}

The Weyl tensor is trace-free and
conformally invariant, ie for $\hat{g}^n=e^{2\phi}g^n$:

\begin{equation}
\label{transweyl}
{W^{\hat{g}^n}_{ijkl}}=e^{2\phi(x)}{W^{g^n}_{ijkl}}
\end{equation}

  While the Schouten tensor has the following transformation law:
\begin{equation}
\label{wschtrans}
P^{\hat{g}^n}_{\alpha\beta}=P^{g^n}_{\alpha\beta}-{\phi}_{\alpha\beta}+
{\phi}_{\alpha}
{\phi}_{\beta}-\frac{1}{2}{\phi}^k{\phi}_kg^n_{\alpha\beta}
\end{equation}

\par In view of our assumption for Theorem \ref{weakt2} and
equation (\ref{Weyl}),
 we may now write $P(g^n)$ in the form:

\begin{equation}
\label{omorfh} P(g^n)=\Sigma_{l\in L} a_l C^l(g^n)
\end{equation}
where each complete contraction $C^l(g^n)$ is in the from:

\begin{equation}
\label{noderis}
contr(W_{i_1j_1k_1l_1}\otimes\dots
\otimes W_{i_Aj_Ak_Al_A}\otimes P_{a_1b_1}
\otimes\dots\otimes P_{a_B b_B})
\end{equation}

 Because of the weight restriction, we see that
$A+B=\frac{n}{2}$.
\newline

\par Let us break up the index set $L$ into subsets
$L^{\mu,\nu}$ as follows: $l\in L^{\mu,\nu}$ if and only if
 $C^l(g^n)$ is
 in the above form and $A=\mu, B=\nu$.

\par We then notice that the linear combination:

$$P^1(g^n)={\Sigma}_{l\in L^{\frac{n}{2},0}} a_l C^l(g^n)$$
is a scalar conformal invariant of weight $-n$. Hence, in view of
the claim of our Theorem \ref{weakt2}, we may subtract it
off, and we are left with considering the case where $P(g^n)$ is a
linear combination:

$$P(g^n)={\Sigma}_{l\in L} a_l C^l(g^n)$$
where each complete contraction $C^l(g^n)$ is in the form
(\ref{noderis}) with $B\ge 1$.

We then have the main theorem of this paper:

\begin{theorem}
\label{Pfaff}
Suppose we are given a $P(g^n)$ which is a linear combination
 of complete contractions of weight $-n$, each in the form
(\ref{noderis}) with $B\ge 1$ and $P(g^n)$ satisfies
(\ref{elie}). Suppose we know
 the coefficient of the complete contraction
$(P^a_a)^{\frac{n}{2}}$ in $P(g^n)$.

\par Then there can be at most one
linear combination $P(g^n)$ of complete contractions in the
 form (\ref{noderis}) with $B\ge 1$ for which the condition
(\ref{elie}) holds.
\end{theorem}

\par If we can show the above, our Theorem
\ref{weakt2} will follow. In order to see this, observe that for
each even dimension $n$, we have that
$\operatorname{Pfaff}(R_{ijkl})$ cannot be
a linear combination of complete contractions depending only
 on the Weyl curvature: If for some $n$ that were the case,
we would
have that for the $n$-sphere $S^n$ with the standard locally
conformally flat metric
$\int_{S^n}\operatorname{Pfaff}(R_{ijkl})dV_{g^n}=0$, which
 is absurd by the Chern-Gauss-Bonnet Theorem.

\par Thus, if we write out $\operatorname{Pfaff}(R_{ijkl})$
 as a linear combination of complete contractions in the form
 (\ref{noderis}) and define
$\overline{\operatorname{Pfaff}}(R_{ijkl})$ to stand for the
sublinear combination of the complete contractions in
$\operatorname{Pfaff}(R_{ijkl})$ with
$B\ge 1$, we will deduce that for some constant $C$,
$P(g^n)$ in Theorem \ref{Pfaff} can be written as:

\begin{equation}
\label{josh}
P(g^n)=C\cdot \overline{\operatorname{Pfaff}}(R_{ijkl})
\end{equation}
This implies our main theorem. $\Box$
\newline

\par We will prove Theorem \ref{Pfaff} by the following two
Lemmas:

\begin{lemma}
\label{A=0}
Given the coefficient of the complete contraction
$(P^a_a)^{\frac{n}{2}}$, there can be at most one sublinear
 combination of complete contractions $C^l(g^n)$ of the form
 (\ref{noderis}) in $P(g^n)$
with $A=0, B=\frac{n}{2}$ so that (\ref{elie}) holds.
\end{lemma}

\begin{lemma}
\label{A>0} Given an integer $1\le A_1\le \frac{n}{2}-1$, and
given the sublinear combination of the complete contractions
$C^l(g^n)$ in $P(g^n)$ with $A< A_1$, then there can be at
 most one sublinear combination of complete contractions
$C^l(g^n)$ of the form (\ref{noderis}) in $P(g^n)$
 with $A=A_1$ so that (\ref{elie}) holds.
\end{lemma}

\par It is clear that if we can prove the above two Lemmas,
 then by induction Theorem \ref{Pfaff} will follow. In the rest of
 the paper we give the proof of these Lemmas.
\newline

\par Our main tool in the proof
will be the super
 divergence formula and the shadow divergence formula used
on the operator $I_{g^n}(\phi)$.

\par A disclaimer on our use of these formulas is in order.
 We will no longer be needing the polarized form
$I^Z_{g^n}(\psi_1,\dots ,\psi_Z)$ of $I^Z_{g^n}(\phi)$. We
 will be referring to the
super divergence formula of $I^Z_{g^n}(\phi)$, and we will mean
the formula that arises from $supdiv[I^Z_{g^n}(\psi_1,\dots
,\psi_Z)]$ by setting $\psi_1=\dots=\psi_Z=\phi$ and dividing by
$Z!$. The same will apply when we refer to the shadow divergence
formula of $I^Z_{g^n}(\phi)$.

\par We must also recall a few more simple facts from
\cite{a:dgciI}. We recall that $I^Z_{g^n}(\phi)$ is taken to be a
linear combination of complete contractions in the form:

\begin{equation}
\label{partlinisym}
\begin{split}
&contr({\nabla}^{m_1}_{r_1\dots r_{m_1}}R_{ijkl}\otimes\dots
\otimes {\nabla}^{m_s}_{t_1\dots t_{m_s}}R_{ijkl}\otimes
\\& {\nabla}^{\nu_1}_{a_1\dots a_{\nu_1}}\phi\otimes
\dots \otimes {\nabla}^{\nu_Z}_{b_1\dots b_{\nu_Z}}\phi)
\end{split}
\end{equation}

\par We also recall that in the context of the iterative
integrations by parts, the $\vec{\xi}$-contractions that we
generically encounter are in the form:

\begin{equation}
 \label{partlinisymxi}
\begin{split}
&contr({\nabla}_{r_1\dots r_{m_1}}^{m_1}R_{i_1j_1k_1l_1}\otimes
\dots \otimes {\nabla}_{v_1\dots
v_{m_s}}^{m_s}R_{i_sj_sk_sl_s}\otimes {\nabla}^{{\nu}_1}_
{{\chi}_1\dots {\chi}_{{\nu}_1}}\phi \otimes\dots \otimes
{\nabla}^{{\nu}_Z}_{{\omega}_1\dots {\omega}_{{\nu}_Z}} \phi\\&
\otimes \vec{\xi}\otimes\dots \otimes \vec{\xi}\otimes
S[{\nabla}^{w_1}_{u_1\dots u_{w_1}}\vec{\xi}]\otimes \dots \otimes
S[{\nabla}^{w_l}_{q_1\dots q_{w_l}}\vec{\xi}] )
\end{split}
\end{equation}
where the factors $\nabla^mR_{ijkl}$ are allowed to have
internal contractions among the indices $i,j,k,l$.

\par Upon occasion, we will be writing those complete
 contractions as linear combinations of complete contractions in the forms:

\begin{equation}
\label{linisym}
\begin{split}
&contr({\nabla}^{m_1}_{r_1\dots r_{m_1}}R_{ijkl}\otimes\dots
\otimes {\nabla}^{m_s}_{t_1\dots t_{m_s}}R_{ijkl}\otimes
\\& S{\nabla}^{\nu_1}_{a_1\dots a_{\nu_1}}\phi\otimes
\dots \otimes S{\nabla}^{\nu_Z}_{b_1\dots b_{\nu_Z}}\phi)
\end{split}
\end{equation}

\begin{equation}
 \label{linisymxi}
\begin{split}
&contr({\nabla}_{r_1\dots r_{m_1}}^{m_1}R_{i_1j_1k_1l_1}\otimes
\dots \otimes {\nabla}_{v_1\dots
v_{m_s}}^{m_s}R_{i_sj_sk_sl_s}\otimes \otimes S{\nabla}^{{\nu}_1}_
{{\chi}_1\dots {\chi}_{{\nu}_1}}\phi \otimes\dots \otimes
\\& S{\nabla}^{{\nu}_Z}_{{\omega}_1\dots
{\omega}_{{\nu}_Z}} \phi
\otimes \vec{\xi}\otimes\dots \otimes \vec{\xi}\otimes
S[{\nabla}^{w_1}_{u_1\dots u_{w_1}}\vec{\xi}]\otimes \dots \otimes
S[{\nabla}^{w_l}_{q_1\dots q_{w_l}}\vec{\xi}] )
\end{split}
\end{equation}
One immediately sees that we can write each complete contraction
in the form (\ref{partlinisym}) or (\ref{partlinisymxi}) as a
linear combination of contractions
 in the forms (\ref{linisym}) or (\ref{linisymxi}) by
repeated use of the identity:

\begin{equation}
\label{curvature}
[\nabla_i\nabla_j-\nabla_j\nabla_i]X_k=R_{ijkl}X^l
\end{equation}

\par We must also recall the
transformation law of the curvature tensor, along with that
 of the Levi-Civita connection,
under conformal re-scalings $\hat{g}^n=e^{2\phi(x)}g^n$:

\begin{equation}
\label{curvtrans}
\begin{split}
&R_{ijkl}^{\hat{g}^n}=e^{2\phi (x)}[R^{g^n}_{ijkl}+ {\phi}_{il}
g_{jk}+{\phi}_{jk}g_{il}-{\phi}_{ik}g_{jl}-{\phi}_{jl}g_{ik}
+{\phi}_i{\phi}_kg_{jl}+{\phi}_j{\phi}_lg_{ik} \\&-{\phi}_i
{\phi}_l g_{jk} -{\phi}_j{\phi}_kg_{il}
+|\nabla\phi|^2g_{il}g_{jk}- |\nabla\phi|^2g_{ik}g_{lj}]
\end{split}
\end{equation}

\begin{equation}
\label{levicivita} {\nabla}^{\hat{g}^n}_k {\eta}_l=
\nabla_k^{g^n}{\eta}_l -{\phi}_k {\eta}_l -{\phi}_l {\eta}_k
+{\phi}^s {\eta}_s g^n_{kl}
\end{equation}

\par Next, we will prove certain Lemmas that will be
  useful throughout this paper.

\subsection{Useful Lemmas.}

Our first Lemma is the following:

\begin{lemma}
\label{firstsymcanc}
 Suppose we are given a collection of
complete contractions $C^k_{g^n}(\phi)$, $k\in K$ of weight
 $-n$ and in the form (\ref{linisym})  or a collection of complete contractions
$C^k_{g^n}(\phi,\vec{\xi})$, $k\in K$, each in the form
(\ref{linisymxi}). Suppose that the identities, respectively:

\begin{equation}
\label{ainte} {\Sigma}_{k\in K} a_k C^k_{g^n}(\phi)=0
\end{equation}

\begin{equation}
\label{ainte2} {\Sigma}_{k\in K} a_k C^k_{g^n}(\phi,
\vec{\xi})=0
\end{equation}

 hold for every Riemannian manifold $(M^n,g^n)$ at any
point $x_0$ and for any function $\phi$ defined around $x_0$, and
in the second case for any vector $\vec{\xi}\in \mathbb{R}^n$.
 We define subsets $K^{(r_1,\dots ,r_Z)}$ of the index set
$K$ as follows:
$k\in K^{(r_1,\dots ,r_Z)}$ if and only if $C^k_{g^n}(\phi)$,
which is in the form (\ref{linisym}), satisfies
${\nu}_1=r_1,\dots ,{\nu}_Z=r_Z$, where the values ${\nu}_1,\dots
,{\nu}_Z$ are taken in decreasing
 rearrangement.

Then, for any subset $K^{(r_1,\dots ,r_Z)}\subset K$, we will
 have, respectively:

\begin{equation}
\label{eqsymcanc} {\Sigma}_{k\in K^{(r_1,\dots ,r_Z)}} a_k
C^k_{g^n}(\phi)=0
\end{equation}

\begin{equation}
\label{eqsymcanc2} {\Sigma}_{k\in K^{(r_1,\dots ,r_Z)}} a_k
C^k_{g^n}(\phi,\vec{\xi}) =0
\end{equation}

 for any Riemannian manifold $(M^n,g^n)$ at any
point $x_0$ and for any function $\phi$ defined around $x_0$, and
in the second case for any vector $\vec{\xi}\in \mathbb{R}^n$.
\end{lemma}

{\it Proof:} We only have to observe that the relations
(\ref{ainte}) and (\ref{ainte2})
 hold formally, where we regard the tensors
$S\nabla^\nu_{r_1\dots r_\nu}\phi$ as symmetric $p$-tensors
$\Omega_{r_1\dots r_\nu}$. On the other hand, the values
${\nu}_1,\dots
,{\nu}_Z$ remain invariant under the permutation relations of
 Definitions 7 and 8 in \cite{a:dgciI}.
 Hence, we have our Lemma. $\Box$

\par Our second Lemma will be the following:

\begin{lemma}
\label{noint} Let us suppose we are given complete contractions
$C^k_{g^n}(\phi)$ in the form (\ref{partlinisym}), and that
 the identity:

$${\Sigma}_{k\in K} a_k C^k_{g^n}(\phi)(x_0)=0$$

holds on any Riemannian manifold $(M^n,g^n)$ and for any
 function  $\phi$ around $x_0$. Let us
suppose that the minimum length among the complete contractions
$\{C^k_{g^n}(\phi)\}_{k\in K}$ is $L$. Then let us define the
subset $K^{\sharp}\subset K$ as follows: $k\in K^{\sharp}$ if
and only if $C^k_{g^n}(\phi)$ which is in the form
(\ref{partlinisym}), has length $L$ and also has
no internal contractions.  We
then have that:

\begin{equation}
\label{nointcanc}
 {\Sigma}_{k\in K^{\sharp}} a_k C^k_{g^n}(\phi)=0
\end{equation}
modulo complete contractions of length $\ge L+1$.
\end{lemma}

{\it Proof:} Let us begin by defining the set $K_1\subset K$
 as
follows: $k\in K_1$ if and only if $C^k_{g^n}(\phi)$ has
 length $L$. Obviously, $K^{\sharp}\subset K_1$.
\par Now, we want to apply Theorem 2 in \cite{a:dgciI}. For each complete
contraction $C^k_{g^n}(\phi),  k\in K_1$, we consider its {\it
linearization} $linC^l(R,\phi)$. Then, by the Lemma hypothesis and
Theorem 2 in \cite{a:dgciI}, we have that the equation:

\begin{equation}
\label{lincanc} {\Sigma}_{k\in K^{\sharp}} a_k linC^k(R,\phi) +
{\Sigma}_{k\in K_1\setminus K^{\sharp}} a_k linC^l(R,\phi)=0
\end{equation}

will hold formally. But then notice the following: For any
linearized complete contraction $linC(R,\phi)$, the number of
internal contractions remains unaltered under any of the
linearized permutation identities. Hence, (\ref{lincanc}) implies
that:

$${\Sigma}_{k\in K^{\sharp}} a_k linC^k(R,\phi)=0$$

formally. But then, as in the proof of the corollaries of
 Theorem 2 in \cite{a:dgciI}, we have that:

$${\Sigma}_{k\in K^{\sharp}} a_k C^k_{g^n}(\phi)=0$$

modulo complete contractions of length $\ge L+1$. $\Box$

\section{The easier step: Proof of Lemma \ref{A=0}.}

Consider any complete contraction $C^l(g^n)$ in the form (\ref{noderis}) with $A=0$.
Let us denote by $R[C^l(g^n)]$ the number of factors $P^a_a$
 in $C^l(g^n)$. Also, let $L^{0,\frac{n}{2},\lambda}$ stand
for the subset of $L$ which is defined as follows:
$l\in L^{0,\frac{n}{2},\lambda}$ if and only if $l\in
L^{0,\frac{n}{2}}$ and $R[C^l(g^n)]=\lambda$.
\newline

\par We will show Lemma \ref{A=0} by an inductive statement.
We assume that for some $T\ge 0$, we have determined the sublinear
combinations $\Sigma_{l\in L^{0,\frac{n}{2},\lambda}} a_l
C^l(g^n)$, for each $\lambda\ge T+1$. We will then show that we
can determine the sublinear combination $\Sigma_{l\in
L^{0,\frac{n}{2},T}} a_l C^l(g^n)$.
 If we can prove this inductive step, then it is obvious
 that our Lemma will follow.
\newline

\par In order to prove the above, we consider
$I^{\frac{n}{2}}_{g^n}(\phi)$.
For any $C^l(g^n)$ with $l\in L^{0,\frac{n}{2}}$, we define
$C^l_{g^n}(\phi)$ to be the complete contraction which is
obtained from $C^l(g^n)$ by substituting each factor $P_{ab}$
 by $-{\nabla}^2_{ab}\phi$.

\par By virtue of (\ref{wschtrans}) and the definition of
$I^{\frac{n}{2}}_{g^n}(\phi)$ we have that:

$$I^{\frac{n}{2}}_{g^n}(\phi)=
{\Sigma}_{l\in L^{0,\frac{n}{2}}} a_l C^l_{g^n}(\phi)$$ modulo
complete contractions of length $\ge \frac{n}{2}+1$. In
particular, each $C^l(g^n)$ with $l\in L^{A,B}$, $A\ge 1$ will not
contribute to the above.
\newline

\par So the problem is reduced to determining the
 sublinear combination
\\ $\Sigma_{l\in L^{0,\frac{n}{2},T}} a_l
C^l_{g^n}(\phi)$ of complete contractions
$C^l_{g^n}(\phi)$
 with $T$ factors $\Delta\phi$ from the sublinear
combination $\Sigma_{s=T+1}^{\frac{n}{2}} \Sigma_{l\in
L^{0,\frac{n}{2},T}} a_l C^l_{g^n}(\phi)$ of complete
contractions $C^l_{g^n}(\phi)$ with more than $T$
 factors $\Delta\phi$.

\par We will use the formula
$supdiv[I^{\frac{n}{2}}_{g^n}(\phi)]$. Let us make a definition:
\newline

Consider any complete contraction $C^l_{g^n}(\phi)$,
$l\in L^{0,\frac{n}{2},T}$. It will be in the form:

$$contr({\nabla}^2_{a_1b_1}\phi\otimes\dots\otimes
{\nabla}^2_{a_{\frac{n}{2}-T}b_{\frac{n}{2}-T}}\phi\otimes
\Delta\phi\otimes\dots\otimes \Delta\phi)$$

where none of the factors ${\nabla}^2_{a_i b_i}\phi$ is in
 the form $\Delta\phi$.

\par We consider the complete contraction $C^{l,D}_{g^n}
(\phi)$:

$$contr({\nabla}^{i_1\dots i_T} [{\nabla}^2_{a_1b_1}\phi\otimes\dots\otimes
{\nabla}^2_{a_{\frac{n}{2}-T}b_{\frac{n}{2}-T}}\phi]\otimes
{\nabla}_{i_1}\phi\otimes\dots \otimes {\nabla}_{i_T}\phi)$$

\par We write out $C^{l,D}_{g^n}(\phi)$
as a linear combination ${\Sigma}_{r\in R^l} a_r
C^r_{g^n}(\phi)$, where each $C^r_{g^n}(\phi)$ is in
the form:

$$contr({\nabla}^{m_1}_{r_1\dots r_{m_1}}
\phi\otimes \dots\otimes {\nabla}^{m_{\frac{n}{2}-T}}_{w_1\dots
w_{m_{\frac{n}{2}-T}}} \phi\otimes{\nabla}_{i_1}\phi
 \otimes\dots \otimes {\nabla}_{i_T}\phi)$$

where each $m_i\ge 2$ and each index $i_s$ contracts against an
index in a factor ${\nabla}^{m_e}\phi$. For each such complete
contraction $C^r_{g^n}(\phi)$, we define $SC^r_{g^n}(\phi)$ to be:

\begin{equation}
\label{symform}
contr(S{\nabla}^{m_1}_{r_1\dots r_{m_1}}
\phi\otimes \dots\otimes S{\nabla}^{m_{\frac{n}{2}-T}}_{w_1\dots
w_{m_{\frac{n}{2}-T}}} \phi\otimes{\nabla}_{i_1}\phi
 \otimes\dots \otimes {\nabla}_{i_T}\phi)
 \end{equation}

\par Observe that, modulo complete contractions of length
$\ge \frac{n}{2}+1$, $C^r_{g^n}(\phi)= SC^r_{g^n}(\phi)$.
\newline

\par For any $l\in L^{0,\frac{n}{2},T}$, we write out
$Tail[C^l_{g^n}(\phi)]$ as a linear combination of complete contractions
in the form (\ref{linisym}). We have that:

\begin{equation}
\label{goodtail}
Tail[C^l_{g^n}(\phi)]= {\Sigma}_{r\in R^l} a_r
SC^r_{g^n}(\phi) +{\Sigma}_{j\in J} a_j C^j_{g^n}(\phi)
\end{equation}

modulo complete contractions of length $\ge \frac{n}{2} +1$.
Each complete contraction $C^j_{g^n}(\phi)$
 has length $\frac{n}{2}$ and less than $T$ factors $\nabla\phi$.

\par Now, for any complete contraction $C^l_{g^n}(\phi)$,
$l\in L^{0,\frac{n}{2},\lambda}$
where $\lambda< T$, we have that:

$$Tail[C^l_{g^n}(\phi)]={\Sigma}_{v\in V} a_v
C^v_{g^n}(\phi)$$
where each complete contraction $C^v_{g^n}(\phi)$
has either length
$\ge \frac{n}{2}+1$ or has length $\frac{n}{2}$ but less than
 $T$ factors $\nabla\phi$. This follows from formula
 (\ref{curvature}).

\par The super divergence formula can be expressed as:

\begin{equation}
\label{supdivana}
\begin{split}
&{\Sigma}_{\lambda =0}^{T-1}{\Sigma}_{l\in L^{0,\frac{n}{2},\lambda}} a_l
Tail[C^l_{g^n}(\phi)]
+ {\Sigma}_{l\in L^{0,\frac{n}{2},T}} a_l
Tail[C^l_{g^n}(\phi)]+
\\& {\Sigma}_{\lambda =T+1}^{\frac{n}{2}-1}
{\Sigma}_{l\in L^{0,\frac{n}{2},\lambda}} a_l
Tail[C^l_{g^n}(\phi)]=0
\end{split}
\end{equation}

modulo complete contractions of length $\ge \frac{n}{2}+1$.
\newline

\par We consider, in (\ref{supdivana}),
the sublinear combination $supdiv[I_{g^n}]|_{\nabla\phi=T}$ of
complete contractions of length $\frac{n}{2}$ with $T$
 factors $\nabla\phi$. From Lemma \ref{firstsymcanc},
we have that

\begin{equation}
\label{bodyguards} supdiv[I_{g^n}]|_{\nabla\phi=T}=0
\end{equation}

\par Furthermore, in view of formula (\ref{supdivana})
and our observations above, we have the following:
Let ${\Sigma}_{\lambda =T+1}^{\frac{n}{2}-1}
{\Sigma}_{l\in L^{0,\frac{n}{2},\lambda}} a_l
Tail[C^l_{g^n}(\phi)]|_{\nabla\phi =T}$
denote the sublinear combination in
${\Sigma}_{\lambda =T+1}^{\frac{n}{2}-1}
{\Sigma}_{l\in L^{0,\frac{n}{2},\lambda}} a_l
Tail[C^l_{g^n}(\phi)]$
of complete contractions with $T$ factors $\nabla\phi$, then:

\begin{equation}
\label{supdivT}
sdI_{\nabla\phi=T}= {\Sigma}_{\lambda =T+1}^{\frac{n}{2}-1}
{\Sigma}_{l\in L^{0,\frac{n}{2},\lambda}} a_l
Tail[C^l_{g^n}(\phi)]|_{\nabla\phi =T}+
 {\Sigma}_{l\in L^{0,\frac{n}{2},T}}a_l
[{\Sigma}_{r\in R^l} a_r
SC^r_{g^n}(\phi)]=0
\end{equation}

 Now, by our inductive hypothesis, we are assuming that we know the
 sublinear combination ${\Sigma}_{\lambda =T+1}^{\frac{n}{2}-1}
{\Sigma}_{l\in L^{0,\frac{n}{2},\lambda}} a_l C^l_{g^n}(\phi)$. Hence,
we deduce that we can determine the sublinear combination
${\Sigma}_{\lambda =T+1}^{\frac{n}{2}-1}
{\Sigma}_{l\in L^{0,\frac{n}{2},\lambda}} a_l Tail[C^l_{g^n}(\phi)]$.
Therefore, we can also determine the sublinear combination
${\Sigma}_{\lambda =T+1}^{\frac{n}{2}-1}
{\Sigma}_{l\in L^{0,\frac{n}{2},\lambda}} a_l Tail[C^l_{g^n}(\phi)]
|_{\nabla\phi=T}$, and using (\ref{supdivT}), we determine the
sublinear combination ${\Sigma}_{l\in L^{0,\frac{n}{2},T}}a_l
[{\Sigma}_{r\in R^l} a_r SC^r_{g^n}(\phi)]$.
\newline

\par A notational convention: When we write $(\nabla)^a$ we
 will mean that we are
taking {\it one} covariant derivative $\nabla_a$ and then raising
the index $a$. (This is to distinguish from $\nabla^a$ which
stands for $a$ iterated covariant derivatives).
 We will now give the following values to factors of
the complete contractions in (\ref{supdivT}):
 To each factor ${\nabla}^2_{ab}\phi$
 we give the value of $-P_{ab}(x_0)$. Also, to each expression of
  the from $S{\nabla}^p_{r_1\dots r_p}\phi ({\nabla})^{r_{i_1}}\phi
\dots ({\nabla})^{r_{i_{p-2}}}\phi$ (where $\{c,d\}=
\{r_1,\dots ,r_p\}\setminus \{r_{i_1},\dots ,r_{i_{p-2}}\}$)
 we give the value $-P_{cd}\cdot (P^a_a)^{p-2}$.
For that assignment $A$ of values, we have that:

$$(n-T)^T\cdot {\Sigma}_{l\in L^{0,\frac{n}{2},T}} a_l C^l(g^n)
 +A\{ {\Sigma}_{\lambda =T+1}^{\frac{n}{2}-1}
{\Sigma}_{l\in L^{0,\frac{n}{2},\lambda}} a_l
Tail[C^l_{g^n}(\phi)]|_{\nabla\phi =T} \}=0$$

\par This concludes the proof of Lemma \ref{A=0}. $\Box$
\newline

\section{The harder step: Proof of Lemma \ref{A>0}.}

\par We want to determine the coefficients of the various
 complete contractions $C^l(g^n)$, indexed in
$L^{A_1,\frac{n}{2}-A_1}$.
\newline

\par We consider $I^{\frac{n}{2}-A_1}_{g^n}(\phi)$. For any
$C^l(g^n)$, $l\in L^{A_1,\frac{n}{2}-A_1}$, we define
$C^l_{g^n}(\phi)$ to be the complete contraction which is obtained
from $C^l(g^n)$ by substituting each factor $P_{ab}$ by
$-{\nabla}_{ab}\phi$.  We then have that:

$$I^{\frac{n}{2}-A_1}_{g^n}(\phi)= {\Sigma}_{l\in L^{A_1,
\frac{n}{2}-A_1}} a_l C^l_{g^n}(\phi)+{\Sigma}_{g\in G} a_g C^g_{g^n}(\phi)$$

modulo complete contractions of length $\ge \frac{n}{2}+1$. The
complete contractions $C^g_{g^n}(\phi)$ are in the form
(\ref{linisym}) and they arise from the sublinear combination
$\Sigma_{k=\frac{n}{2}-A_1 +1}^{\frac{n}{2}} \Sigma_{l\in L^{
\frac{n}{2}-k,k}} a_l C^l(g^n)$. Hence, we have that the sublinear
combination $\Sigma_{g\in G} a_g C^g_{g^n}(\phi)$ is known.
\newline

\par The complete contractions $C^l_{g^n}(\phi)$ are in the form :

\begin{equation}
\label{weylphi}
contr(W_{i_1j_1k_1l_1}\otimes\dots \otimes W_{i_{A_1}j_{A_1}
k_{A_1}l_{A_1}}\otimes {\nabla}^2_{a_1b_1}\phi\otimes\dots
 \otimes {\nabla}^2_{a_{\frac{n}{2}-A_1}
b_{\frac{n}{2}-A_1}}\phi)
\end{equation}

\par While we write the complete contractions $C^g_{g^n}(\phi)$ in the form:

\begin{equation}
\label{diakrita}
contr(R_{i_1j_1k_1l_1}\otimes\dots\otimes R_{i_{A_1}j_{A_1}
k_{A_1}l_{A_1}}\otimes {\nabla}^2_{a_1b_1}\phi\otimes\dots \otimes {\nabla}^2_{a_{\frac{n}{2}-A_1}
b_{\frac{n}{2}-A_1}}\phi)
\end{equation}
(for this equation, the factors $R_{ijkl}$, $\nabla^2_{ab}\phi$
are allowed to have internal contractions).

\par Now, we write ${\Sigma}_{l\in L^{A_1,
\frac{n}{2}-A_1}} a_l C^l_{g^n}(\phi)$
as a linear combination:

\begin{equation}
\label{brent}
{\Sigma}_{l\in L^{A_1,
\frac{n}{2}-A_1}} a_l C^l_{g^n}(\phi)=\Sigma_{u\in U} a_u C^u_{g^n}
(\phi)
\end{equation}
where each $C^u_{g^n}(\phi)$ is in the form:

\begin{equation}
\label{refback}
\begin{split}
&contr(R_{ijkl}\otimes\dots R_{i'j'k'l'}\otimes Ric_{kl}
\otimes \dots \otimes Ric_{k'l'}\otimes R\otimes \dots \otimes R\otimes
\\& {\nabla}^2_{\alpha\beta}\phi\otimes \dots
\otimes{\nabla}^2_{{\alpha}'{\beta}'}\phi \otimes \Delta\phi
\otimes\dots \otimes\Delta\phi)
\end{split}
\end{equation}
When we employ the above notation we will imply that
 each of the factors $R_{ijkl}$, $Ric_{ab}$ and
${\nabla}^2_{\alpha\beta}\phi$  does not have any of the indices
 $i,j,k,l$ or $a,b$ or
$\alpha,\beta$ contracting between themselves. Let $Z$ stand for
the number of factors $R_{ijkl}$, $X$ for the number of factors
$Ric_{ab}$,
 $C$ for the number of factors $R$, $\Gamma$ for the number
 of factors ${\nabla}^2_{\alpha\beta}\phi$ and $\Delta$
for the number of factors $\Delta\phi$. We have that
$Z+X+C=A_1$ and $\Gamma+\Delta =\frac{n}{2}-A_1$.

\par We denote the corresponding index set in $U$ by
$U^{Z,X,C,\Gamma,\Delta}$. We then claim the following:

\begin{lemma}
\label{chuck}
\par Under the assumptions of Lemma \ref{A>0}, we claim that we can
determine all the sublinear combinations
$\Sigma_{u\in U^{Z,X,C,\Gamma,\Delta}} a_u C^u_{g^n}(\phi)$ above.
\end{lemma}

\par Before we prove this Lemma, let us explain how we can deduce
our desired Lemma \ref{A>0} from Lemma \ref{chuck}.

\par If we can determine all the sublinear combinations
$\Sigma_{u\in U^{Z,X,C,\Gamma,\Delta}} a_u C^u_{g^n}(\phi)$, we
then will have determined the whole linear combination
$\Sigma_{u\in U} a_u C^u_{g^n}(\phi)$, and hence by (\ref{brent})
we will have determined the linear combination ${\Sigma}_{l\in
L^{A_1,\frac{n}{2}-A_1}} a_l C^l_{g^n}(\phi)$.

\par But then, setting $\nabla^2_{ab}\phi (x_0)=-P_{ab}(x_0)$, we determine
 ${\Sigma}_{l\in L^{A_1,
\frac{n}{2}-A_1}} a_l C^l(g^n)$, and hence we will have shown our
 Lemma.
\newline

\subsection{The long induction: The Proof of Lemma \ref{chuck}.}

\par We will determine the various sublinear combinations by an
induction.
\newline

\par We initially determine the  sublinear combination
$\Sigma_{u\in U^{0,1,A_1-1,1,\frac{n}{2}-A_1-1}} a_u C^u_{g^n}(\phi)$.
 By definition, we see that the sublinear combination in question will
 be of the form $(const)\cdot C^{*}_{g^n}(\phi)$, where
$C^{*}_{g^n}(\phi)$ is the complete contraction:

\begin{equation}
\label{hbre}
contr(R^{A_1-1}\otimes Ric^{ab}\otimes \nabla^2_{ab}\phi\otimes
(\Delta\phi)^{\frac{n}{2}-A_1-1})
\end{equation}
(Thus, determining $\Sigma_{u\in
U^{0,1,A_1-1,1,\frac{n}{2}-A_1-1}} a_u C^u_{g^n}(\phi)$ amounts to
determining $(const)$).
\newline

\par Then, we will determine the sublinear combination
$\Sigma_{u\in U^{0,0,A_1, 0,\frac{n}{2}-A_1}} a_u C^u_{g^n}(\phi)$. We
 observe that this sublinear combination will be in the form:

\begin{equation}
\label{hbre2}
(const)'\cdot contr(R^{A_1}\otimes (\Delta\phi)^{\frac{n}{2}-A_1})
\end{equation}
(Thus again, we only have to determine $(const)'$).
\newline

\par Finally, having determined the two sublinear combinations above, we
 will prove the following inductive statement: Let us suppose that for
some number $\Delta_1+1$, we have determined all the sublinear
 combinations
$\Sigma_{u\in U^{Z,X,C,\Gamma,\Delta}} a_u C^u_{g^n}(\phi)$
with $\Delta\ge \Delta_1+1$. Moreover, we assume that for some number
 $C_1+1$, we have determined all the sublinear combinations
$\Sigma_{u\in U^{Z,X,C,\frac{n}{2}-A_1-\Delta_1,\Delta_1}} a_u
C^u_{g^n}(\phi)$ with $C\ge C_1+1$. Finally, we
 suppose that for some number $X_1+1$, we have determined all the
 sublinear combinations $\Sigma_{u\in U^{Z,X,C_1,
\frac{n}{2}-A_1-\Delta_1,\Delta_1}} a_u C^u_{g^n}(\phi)$ with
$X\ge X_1+1$. We then claim that we can determine the sublinear
 combination $\Sigma_{u\in U^{A_1-X_1-C_1,X_1,C_1,
\frac{n}{2}-A_1-\Delta_1,\Delta_1}} a_u C^u_{g^n}(\phi)$. If we can
 show  the above then by induction we will have proven our Lemma
\ref{chuck}.
\newline

\par Before proceeding with the proof, we make note of how the
Weyl tensor can be decomposed:

\begin{equation}
\label{weyl''}
\begin{split}
&W_{ijkl}=R_{ijkl}+ \frac{1}{n-2}[Ric_{ik}g^n_{jl}+Ric_{jl}g^n_{ik}
-Ric_{il}g^n_{jk}-Ric_{jk}g^n_{il}]
\\& -\frac{R}{(n-1)(n-2)}g^n_{ik}
g^n_{jl} +\frac{R}{(n-1)(n-2)}g^n_{il}g^n_{jk}
\end{split}
\end{equation}

{\it Determining the sublinear combination
$\Sigma_{u\in U^{0,1,A_1-1,1,\frac{n}{2}-A_1-1}} a_u C^u_{g^n}(\phi)$:}
\newline

\par  We consider
$I^{\frac{n}{2}-A_1+1}_{g^n}(\phi)$. We focus on the sublinear
 combinations of complete contractions of length $\frac{n}{2}$ or
$\frac{n}{2} +1$
 in $I^{\frac{n}{2}-A_1+1}_{g^n}(\phi)$, which we respectively denote
 by $I^{\frac{n}{2}-A_1+1}_{g^n}(\phi)|_{\frac{n}{2}}$,
$I^{\frac{n}{2}-A_1+1}_{g^n}(\phi)|_{\frac{n}{2} +1}$. Using the
transformation law (\ref{wschtrans}) and the conformal invariance of
 the Weyl tensor, we deduce that the sublinear combination
$I^{\frac{n}{2}-A_1+1}_{g^n}(\phi)|_{\frac{n}{2}}$ arises from the
 sublinear combination $\Sigma_{B=0}^{A_1-1}
\Sigma_{l\in L^{B,\frac{n}{2}-B}}a_l C^l(g^n)$ in
$P(g^n)$. Therefore by our inductive hypothesis, we have that the
 sublinear combination $I^{\frac{n}{2}-A_1+1}_{g^n}(\phi)
|_{\frac{n}{2}}$ in
$I^{\frac{n}{2}-A_1+1}_{g^n}(\phi)$ is known.
\newline

\par Now, we also claim that the sublinear combination
$I^{\frac{n}{2}-A_1+1}_{g^n}(\phi)|_{\frac{n}{2} +1}$ in
\\$I^{\frac{n}{2}-A_1+1}_{g^n}(\phi)$ can be written as:

\begin{equation}
\label{pfafanalys}
I^{\frac{n}{2}-A_1+1}_{g^n}(\phi)|_{\frac{n}{2} +1}=\Sigma_{k\in K} a_k
C^k_{g^n}(\phi) +\Sigma_{u\in U_1} a_u C^u_{g^n}(\phi)
\end{equation}
where  $\Sigma_{k\in K} a_k
C^k_{g^n}(\phi)$ arises from
$\Sigma_{B=0}^{A_1-1} \Sigma_{l\in L^{B,\frac{n}{2}-B}} a_l
C^l(g^n)$ in $P(g^n)$ and $\Sigma_{u\in U_1} a_u C^u_{g^n}(\phi)$ arises from the sublinear
combination $\Sigma_{l\in L^{A_1,\frac{n}{2}-A_1}}a_l C^l(g^n)$ in
$P(g^n)$. This means that the contractions $C^l(g^n)$, $l\in L^{B,
\frac{n}{2}-B}$ with $B\ge A_1+1$ will not contribute to
$I^{\frac{n}{2}-A_1+1}_{g^n}(\phi)|_{\frac{n}{2}+1}$. This follows
by virtue of (\ref{wschtrans}). Hence, we may assume that the
sublinear combination $\Sigma_{k\in K} a_k C^k_{g^n}(\phi)$ is
known.

\par Now, we initially have that the complete contractions
$C^u_{g^n}(\phi)$ on the right hand side of the above are in the form:

\begin{equation}
\label{curvphi2}
contr(W_{i_1j_1k_1l_1}\otimes\dots \otimes W_{i_{A_1}j_{A_1}
k_{A_1}l_{A_1}}\otimes {\nabla}^2_{a_1b_1}\phi\otimes\dots
 \otimes {\nabla}^2_{a_{\frac{n}{2}-A_1-1}
b_{\frac{n}{2}-A_1-1}}\phi\otimes \nabla_x\phi\otimes \nabla_d\phi)
\end{equation}

\par Then, we decompose the Weyl tensor as in (\ref{weyl''}) and we
write the linear combination on the right hand side of the above as a
linear combination of complete contractions in the form:

\begin{equation}
\label{kakomoirh}
\begin{split}
&contr(R_{i_1j_1k_1l_1}\otimes\dots \otimes R_{i_zj_z
k_zl_z}\otimes Ric_{h_1e_1}\otimes\dots\otimes Ric_{h_ye_y}\otimes
R^q\otimes {\nabla}^2_{a_1b_1}\phi\otimes\dots
\\& \otimes {\nabla}^2_{a_qb_q}\phi\otimes
 \nabla_x\phi\otimes\nabla_d\phi \otimes (\Delta\phi)^r)
\end{split}
\end{equation}
where we are making the notational convention that no
two indices in any factor $R_{ijkl}, Ric_{ab}, \nabla^2_{ij}\phi$ are
contracting between themselves. We write:

$$\Sigma_{u\in U_1} a_u C^u_{g^n}(\phi)=
\Sigma_{u\in U_2} a_u C^u_{g^n}(\phi)$$
where each $C^u_{g^n}(\phi)$, $u\in U_2$ is in the form
(\ref{kakomoirh}). We replace the expression
$\Sigma_{u\in U_2} a_u C^u_{g^n}(\phi)$ for
$\Sigma_{u\in U_1} a_u C^u_{g^n}(\phi)$ in (\ref{pfafanalys}).
 Moreover, we assume that each $C^k_{g^n}(\phi)$ in (\ref{pfafanalys})
 is in the form (\ref{kakomoirh}).

\par Now, we focus on the sublinear combination in
$\Sigma_{u\in U_2} a_u C^u_{g^n}(\phi)$ that consists of complete
contractions in the form (\ref{kakomoirh}) with $Z=0$ factors
$R_{ijkl}$, $Y=1$ factor $Ric_{he}$, $C=A_1-1$ factors $R$,
$\Gamma=0$ factors $\nabla^2\phi$, $\Delta=\frac{n}{2}-A_1 -1$
 factors $\Delta\phi$. We also assume that the two factors $\nabla\phi$
 contract against the two indices of the one factor $Ric_{ij}$.
 Therefore, we
 have that the sublinear combination in question is of the form
$(const)_{*}\cdot C^{*}_{g^n}(\phi)$, where $C^{*}_{g^n}(\phi)$ is in
 the form:

\begin{equation}
\label{starkri9ar}
contr(R^{A_1-1}\otimes Ric^{ij}\otimes\nabla_i\phi\otimes\nabla_j\phi
\otimes (\Delta\phi)^{\frac{n}{2}-A_1-1})
\end{equation}

\par We now make two claims:

\begin{lemma}
\label{enadyo}
We have that the sublinear combination $(const)_{*}\cdot
C^{*}_{g^n}(\phi)$ in $\Sigma_{u\in U_2} a_u C^u_{g^n}(\phi)$
 arises from the sublinear combination
$\Sigma_{u\in U^{0,1,A_1-1,1,\frac{n}{2}-A_1-1}} a_u C^u_{g^n}(\phi)$
in $I^{\frac{n}{2}-A_1}_{g^n}(\phi)$ by replacing the factor
$\nabla^2_{ij}\phi$ by an expression $-\nabla_i\phi\nabla_j\phi$.

\par Our second claim is that the sublinear combination
$(const)_{*}\cdot C^{*}_{g^n}(\phi)$ can be determined from the
known sublinear combinations in (\ref{pfafanalys}), using the
shadow divergence formula for
$I^{\frac{n}{2}-A_1+1}_{g^n}(\phi)$.
\end{lemma}

\par We observe that if we can show the above Lemma, we will then have determined the sublinear combination
$\Sigma_{u\in U^{0,1,A_1-1,1,\frac{n}{2}-A_1-1}} a_u
C^u_{g^n}(\phi)$ in $I^{\frac{n}{2}-A_1}_{g^n}(\phi)$, and hence
proven the first base case of our induction.
\newline

{\it Proof of Lemma \ref{enadyo}:} We begin with the first part.
Initially, let us focus on the sublinear combination
$\Sigma_{u\in U^{0,1,A_1-1,1,\frac{n}{2}-A_1-1}} a_u C^u_{g^n}(\phi)$ in
$I^{\frac{n}{2}-A_1}_{g^n}(\phi)$ and understand in detail how it
 arises. For each $l\in L^{A_1,\frac{n}{2}-A_1}$, we consider the
 complete contraction $C^l_{g^n}(\phi)$ defined above, which will be in
 the form (\ref{weylphi}). We then decompose the factors $W_{ijkl}$ as
 in (\ref{weyl''}).

\par Now, for each factor $W_{ijkl}$, we have the
 option of replacing it by one of the 7 expressions on the right hand
 side of (\ref{weyl''}). Therefore, we can write $C^l_{g^n}(\phi)$ as a
sum of $7^{A_1}$ complete contractions in the form (\ref{diakrita}):

\begin{equation}
\label{mihelich} C^l_{g^n}(\phi)=\Sigma_{\tau =1}^{7^{A_1}} a_\tau
C^\tau_{g^n}(\phi)
\end{equation}

 Each of the $7^{A_1}$ different summands corresponds to a different
 sequence of substitutions of the $A_1$ factors $W_{ijkl}$ as explained
 above. We then group up the complete contractions $C^\tau_{g^n}(\phi)$
on the right hand side of the above that are of the form
(\ref{hbre}), and we
 denote that sublinear combination in (\ref{mihelich}) by
$F[C^l_{g^n}(\phi)]$. Hence, using this notation we have that:

$$\Sigma_{u\in U^{0,1,A_1-1,1,\frac{n}{2}-A_1-1}} a_u C^u_{g^n}(\phi)=
\Sigma_{l\in L^{A_1,\frac{n}{2}-A_1}} a_l F[C^l_{g^n}(\phi)]$$

\par Now, we consider the complete contractions in
$Image^{\frac{n}{2}-A_1+1}_\phi[C^l(g^n)]$, for each
$l\in L^{A_1,\frac{n}{2}-A_1}$. We are only interested in the sublinear combination
$$Image^{\frac{n}{2}-A_1+1}_\phi|_{\frac{n}{2} +1}[C^l(g^n)]$$

of complete contractions of length $\frac{n}{2} +1$. It follows
that this sublinear combination arises by replacing $\frac{n}{2}
-A_1-1$
 factors $P_{ab}$ by the  expression
$-\nabla^2_{ab}\phi$ on the right hand side of (\ref{wschtrans})
and
also by replacing one factor $P_{ab}$ by a quadratic expression on the
right hand side of (\ref{wschtrans}).

\par Now, we further denote by
$Image^{\frac{n}{2}-A_1+1,+}_\phi|_{\frac{n}{2}+1}[C^l(g^n)]$
the sublinear combination in
$Image^{\frac{n}{2}-A_1+1}_\phi|_{\frac{n}{2}+1}[C^l(g^n)]$ that arises
when we replace $\frac{n}{2} -A_1-1$ factors $P_{ab}$ by
$-\nabla^2_{ab}\phi$ and one factor $P_{ab}$ by the expression
$g_{ab}|\nabla\phi|^2$. We trivially observe that if we write out
$Image^{\frac{n}{2}-A_1+1,\sigma +1,+}_\phi[C^l(g^n)]$
as a linear combination of complete contractions in the form
(\ref{kakomoirh}), none will be in the form (\ref{starkri9ar}).

\par Hence, we may restrict our attention to the sublinear combination
\\ $Image^{\frac{n}{2}-A_1+1,-}_\phi[C^l(g^n)]$
 in $Image^{\frac{n}{2}-A_1+1}_\phi[C^l(g^n)]$
that arises when we replace $\frac{n}{2}-A_1-1$ factors $P_{ab}$ by
 $-\nabla^2_{ab}\phi$ and one factor $P_{ab}$ by $\nabla_a\phi
\nabla_b\phi$. Hence, comparing $Image^{\frac{n}{2}-A_1}\{
\Sigma_{l\in L^{A_1,\frac{n}{2}-A_1}} a_l C^l(g^n)\}$ and
$Image^{\frac{n}{2}-A_1+1,-}\{
\Sigma_{l\in L^{A_1,\frac{n}{2}-A_1}} a_l C^l(g^n)\}$,
we see that $Image^{\frac{n}{2}-A_1+1,-}\{
\Sigma_{l\in L^{A_1,\frac{n}{2}-A_1}} a_l C^l(g^n)\}$ arises from
$Image^{\frac{n}{2}-A_1}\{
\\ \Sigma_{l\in L^{A_1,\frac{n}{2}-A_1}} a_l C^l(g^n)\}$
by picking out
one factor $\nabla^2_{ab}\phi$ from each complete contraction in the
form (\ref{weylphi}) in  $Image^{\frac{n}{2}-A_1}\{
\Sigma_{l\in L^{A_1,\frac{n}{2}-A_1}} a_l C^l(g^n)\}$
 (this factor may now also be of the form $\Delta\phi$)
  and replacing it by an expression $-\nabla_a\phi\nabla_b\phi$.
In that case, if we repeat the decomposition of the
 factors $W_{ijkl}$ to the complete contractions in
 $Image^{\frac{n}{2}-A_1+1,-}_\phi[C^l(g^n)]$, we obtain the
first claim of our Lemma.
\newline

\par Now, for the second part of our Lemma, we first of all denote
$(const)_{*} C^{*}_{g^n}(\phi)$ by $\Sigma_{u\in U^{*}_2} a_u
C^u_{g^n}(\phi)$. We then want to apply the shadow divergence formula
 to $I^{\frac{n}{2}-A_1+1}_{g^n}(\phi)$ and determine
the sublinear combination $\Sigma_{u\in U^{*}_2} a_u C^u_{g^n}(\phi)$.
 We will focus on the sublinear combination of $\vec{\xi}$-contractions
 in $Shad[I^{\frac{n}{2}-A_1+1}_{g^n}(\phi)]$ that are in the form:

\begin{equation}
\label{hellenicus}
contr((|\vec{\xi}|^2)^{A_1-1}\otimes
S\nabla^{\frac{n}{2}-A_1}_{r_1\dots r_{\frac{n}{2}-A_1}}
\vec{\xi}_{r_{\frac{n}{2}-A_1+1}}\otimes\nabla^{r_1}\phi\otimes\dots
\otimes \nabla^{r_{\frac{n}{2}-A_1+1}}\phi)
\end{equation}

\par If we denote the sublinear combination of those
$\vec{\xi}$-contractions in $Shad[I^{\frac{n}{2}-A_1+1}_{g^n}(\phi)]$
by $Shad_{+}[I^{\frac{n}{2}-A_1+1}_{g^n}(\phi)]$, we claim that:

\begin{equation}
\label{brent4}
Shad_{+}[I^{\frac{n}{2}-A_1+1}_{g^n}(\phi)]=0
\end{equation}

\par This is straightforward because the shadow divergence formula holds
 formally. Now, for each $k\in K$ (see (\ref{pfafanalys})) we denote
 by $Tail^{Shad}_{+}[C^k_{g^n}(\phi)]$ the sublinear combination in each
 $Tail^{Shad}[C^k_{g^n}(\phi)]$ that consists of
$\vec{\xi}$-contractions in the form (\ref{hellenicus}). Analogously, for each $u\in U_2$, we denote by
$Tail^{Shad}_{+}[C^u_{g^n}(\phi)]$ the sublinear combination in each
 $Tail^{Shad}[C^u_{g^n}(\phi)]$ that consists of
$\vec{\xi}$-contractions in the form (\ref{brent4}).
Now, we observe that the $\vec{\xi}$-length of
 the $\vec{\xi}$-contraction in (\ref{brent4}) is $\frac{n}{2} +1$.
 Hence, in view of the Lemma on acceptable descendants in
\cite{a:dgciI} and also (\ref{pfafanalys}), (\ref{brent4}),
 we deduce that:

\begin{equation}
\label{brent5}
Shad_{+}[I^{\frac{n}{2}-A_1+1}_{g^n}(\phi)|_{\frac{n}{2}}]
=Tail^{Shad}_{+}[\Sigma_{k\in K} a_k C^k_{g^n}(\phi)]
+Tail^{Shad}_{+}[\Sigma_{u\in U_2} a_u C^u_{g^n}(\phi)]=0
\end{equation}

\par Therefore, if we could show that for each $u\in U_2\setminus
U^{*}_2$, we have that:

\begin{equation}
\label{brent6}
 Shad_{+}[C^u_{g^n}(\phi)]=0
\end{equation}
we could then use equation (\ref{brent5}) to determine the sublinear
 combination
$$\Sigma_{u\in U^{*}_2} a_u Tail_{+}^{Shad}
[C^u_{g^n}(\phi)].$$

\par Let us observe how
it would then be straightforward to determine
 $\Sigma_{u\in U^{*}_2} a_u C^u_{g^n}(\phi)$:

\par We claim that for $u\in U^{*}_2$,
$Tail_{+}^{Shad}[C^u_{g^n}(\phi)]=(-1)^{\frac{n}{2}-A_1}C^{*}_{g^n}
(\phi,\vec{\xi})$. To see this, we note that
$Tail_{+}^{Shad}[C^u_{g^n}(\phi)]$ arises in the following way:
Let us denote by $C^{*}_{g^n}(\psi,\vec{\xi})$ the descendant of
$C^u_{g^n}(\phi)$ that arises by replacing the $A_1-1$ factors $R$
by $|\vec{\xi}|^2$, the $\frac{n}{2}-A_1-1$ factors
$\Delta\phi$
by $\nabla_i\phi \vec{\xi}^i$ and the factor $Ric_{ij}$ by
$-\nabla_i\vec{\xi}_j$ (all in $N$-cancelled notation). Recall
from \cite{a:dgciI} that $O^{Shad}[C_{g^n}(\phi,\vec{\xi})]$ stands for the
 sublinear combination of hard and stigmatized $\vec{\xi}$-contractions
(of both types) that arise along the iterative
integrations by parts of the $\vec{\xi}$-contraction
$C_{g^n}(\phi,\vec{\xi})$. We now define $O^{Shad}_{+}
[C_{g^n}(\phi,\vec{\xi})]$ to stand for the sublinear combination of
 those $\vec{\xi}$-contractions that are in the form
(\ref{hellenicus}).

We then claim (claim 1) that $O^{Shad}_{+}
[C^{*}_{g^n}(\phi,\vec{\xi})]=(-1)^{\frac{n}{2}-A_1}C^{*}_{g^n}
(\phi,\vec{\xi})$ (the left hand side here stands for the sublinear
 combination of complete contractions in the form
$C^{*}_{g^n}(\phi,\vec{\xi})$ in $O^{Shad}
[C^{*}_{g^n}(\phi,\vec{\xi})]$).
 Moreover, we claim (claim 2) that for any other descendant
$C^d_{g^n}(\phi,\vec{\xi})$ of $C^{*}_{g^n}(\phi)$ we will have
$O^{Shad}_{+}[C^d_{g^n}(\phi,\vec{\xi})]=0$.

\par The second claim follows
by simply observing that $C^d_{g^n}(\phi,\vec{\xi})$ must contain a
 factor with an internal contraction, hence each $\vec{\xi}$-contraction
 in $O^{Shad}[C^d_{g^n}(\phi,\vec{\xi})]$ with length
$\frac{n}{2} +1$ must have a factor with an internal contraction.
Our first claim follows by integrating by parts all the factors
$\vec{\xi}_i$ that contract against factors $\nabla\phi$ and making
all the derivatives $\nabla_i$ hit the factor $\nabla\vec{\xi}$ and then
 symmetrizing. We observe that any other $\vec{\xi}$-contraction that
arises in the iterative integration by parts will not be of the
form $C^{*}_{g^n}(\phi,\vec{\xi})$: It will either have
$\vec{\xi}$-length $\ge\frac{n}{2}+2$ or a factor $\nabla^a\phi,
a\ge 2$ or less than $A_1-1$ factors $|\vec{\xi}|^2$.

\par In view of the above, and since
 (\ref{brent5}) holds formally, if we replace each expression
$S\nabla^{\frac{n}{2}-A_1}_{r_1\dots r_{\frac{n}{2}-A_1}}
\vec{\xi}_{r_{\frac{n}{2}-A_1+1}}\otimes\nabla^{r_1}\phi\otimes\dots
\otimes \nabla^{r_{\frac{n}{2}-A_1+1}}\phi$
 in each complete contraction in (\ref{brent5}) by
$\nabla^i\phi\nabla^j\phi\nabla_i\vec{\xi}_j
(\Delta\phi)^{\frac{n}{2}-A_1-1}$ and each factor $|\vec{\xi}|^2$ by a
 factor $R$, we can then determine $\Sigma_{u\in U^{*}_2} a_u
C^u_{g^n}(\phi)$.
Hence, showing (\ref{brent6}) would complete the proof of our
Lemma.
\newline

\par But (\ref{brent6}) is easy to prove: Let us suppose that
$C^u_{g^n}(\phi)$ is in the form (\ref{kakomoirh}) and
has less than $A_1-1$ factors $R$. It then follows that each descendent
$C^{u,l}_{g^n}(\phi,\vec{\xi})$ of $C^u_{g^n}(\phi)$ will have less
 than $A_1-1$ factors $|\vec{\xi}|^2$ (by the Lemma on the acceptable descendants
 in \cite{a:dgciI}) and
 hence, by the iterative integration by parts procedure, each
 $\vec{\xi}$-contraction in $Tail^{Shad}[C^{u,l}_{g^n}
(\phi,\vec{\xi})]$ will have less than $A_1-1$ factors
$|\vec{\xi}|^2$ (by Lemma 15 in \cite{a:dgciI}) and hence we have
shown (\ref{brent6}) in this
 case. Now, we consider the case where $C^u_{g^n}(\phi)$ is in the form
(\ref{kakomoirh}) and has less than $\frac{n}{2}-A_1-1$ factors
$\Delta\phi$, and hence has at least one factor
$\nabla^2\phi\ne \Delta\phi$.
 It then follows that each descendent
$C^{u,l}_{g^n}(\phi,\vec{\xi})$ of $C^u_{g^n}(\phi)$
will have at least one factor $\nabla^2\phi\ne \Delta\phi$
 (by the Lemma on the acceptable descendants
 in \cite{a:dgciI}). Hence,
we have that each $\vec{\xi}$-contraction of $\vec{\xi}$-length
$\frac{n}{2} +1$ in $Tail^{Shad}[C^u_{g^n}(\phi)]$ will have at
least
 one factor $\nabla^a\phi, a\ge 2$ and therefore
$Tail^{Shad}_{+}[C^u_{g^n}(\phi)]=0$.

\par We are thus left with the case where $u\in U_2\setminus U^{*}_2$
 and $C^u_{g^n}(\phi)$ has at least $A_1-1$ factors $R$ and at least
$\frac{n}{2}-A_1-1$ factors $\Delta\phi$. It then follows that
$C^u_{g^n}(\phi)$ must be in the form:

$$contr (R^{A_1}\otimes (\Delta\phi)^{\frac{n}{2}-A_1-1}\otimes
|\nabla\phi|^2)$$
 But then, by the iterative integrations by parts procedure, we observe
 that each $\vec{\xi}$-contraction of $\vec{\xi}$-length
$\frac{n}{2} +1$ in $Tail^{Shad}[C^u_{g^n}(\phi)]$ will either have a
 factor $\nabla^a\phi, a\ge 2$ or will have two factors
$\nabla\phi$
 that contract against each other. Therefore, we again have our desired
(\ref{brent6}) in this case. We have shown our Lemma. $\Box$
\newline

{\it Determining the sublinear combination
$\Sigma_{u\in U^{0,0,A_1,0,\frac{n}{2}-A_1}} a_u C^u_{g^n}(\phi)$:}
\newline

We consider the shadow divergence formula of
$I^{\frac{n}{2}-A_1}_{g^n}(\phi)$,
$Shad[I^{\frac{n}{2}-A_1}_{g^n}(\phi)]$.
We focus on the sublinear combination of $\vec{\xi}$-contractions in
 the form:

\begin{equation}
\label{stewart}
contr((|\vec{\xi}|^2)^{A_1-1}\otimes
\nabla^{\frac{n}{2}-A_1+1}_{r_1\dots
r_{\frac{n}{2}-A_1-1}i}\vec{\xi}_j\otimes \nabla^{ij}\phi\otimes
(\nabla)^{r_1}\phi\otimes\dots \otimes
(\nabla)^{r_{\frac{n}{2}-A_1-1}}\phi)
\end{equation}

 We denote the above $\vec{\xi}$-contraction by
$C^\sharp_{g^n}(\phi,\vec{\xi})$ for short.
For each $C_{g^n}(\phi)$ in the form (\ref{diakrita}) of length
$\frac{n}{2}$, we denote
 by $Tail^{Shad}_{+}[C_{g^n}(\phi)]$ the sublinear combination of
$\vec{\xi}$-contractions in the form (\ref{stewart}) in $Tail^{Shad}
[C_{g^n}(\phi)]$. (Note that we are changing the meaning of
$Tail^{Shad}_{+}[C_{g^n}(\phi)]$).
This notation extends to linear combinations.
Now, since the Shadow divergence formula holds formally, we will have
 that:

$$Tail^{Shad}_{+}[I^{\frac{n}{2}-A_1}_{g^n}(\phi)]=0$$

\par We write out $I^{\frac{n}{2}-A_1}_{g^n}(\phi)$ in the
 form:

\begin{equation}
\label{baskin}
I^{\frac{n}{2}-A_1}_{g^n}(\phi)=\Sigma_{k\in K} a_k C^k_{g^n}(\phi)+
\Sigma_{u\in U} a_u C^u_{g^n}(\phi)
\end{equation}
modulo complete contractions of length $\ge \frac{n}{2}+1$. Here $\Sigma_{k\in K} a_k C^k_{g^n}(\phi)$
arises from the sublinear combination
$\Sigma_{A=0}^{A_1-1}\Sigma_{l\in L^{A,\frac{n}{2}-A}} a_l
C^l(g^n)$. Hence, we have that $\Sigma_{k\in K} a_k
C^k_{g^n}(\phi)$ is known. We note that the index set $K$ differs
from $K$ in (\ref{pfafanalys}). We deduce that:

\begin{equation}
\label{stuart2}
\Sigma_{k\in K} a_k Tail^{Shad}_{+} [C^k_{g^n}(\phi)]+
\Sigma_{u\in U} a_u Tail^{Shad}_{+} [C^u_{g^n}(\phi)]=0
\end{equation}

\par Now, we claim that for each $u\in  U\setminus (U^{0,0,A_1,0,
\frac{n}{2}-A_1}\bigcup U^{0,1,A_1-1,1,\frac{n}{2}-A_1 -1})$ we have
that $Tail^{Shad}_{+} [C^u_{g^n}(\phi)]=0$. This follows by a similar
reasoning as for the previous case: For each $u$ above, we have that
 either $C^u_{g^n}(\phi)$ has less than $A_1-1$ factors $R$ or it has
 less than $\frac{n}{2}-A_1-1$ factors $\Delta\phi$.
In the first case we
 then have that each $\vec{\xi}$-contraction in
$Tail^{Shad}[C^u_{g^n}(\phi)]$ will have less than
$A_1-1$ factors $|\vec{\xi}|^2$ and in the second, it will have less
 than $\frac{n}{2}-A_1-1$ factors $\nabla\phi$.

\par In view of this fact, we can then use
(\ref{stuart2}) to determine the sublinear combination
$\Sigma_{u\in U^{0,0,A_1,0,\frac{n}{2}-A_1}} a_u
Tail^{Shad}_{+} [C^u_{g^n}(\phi)]$.
\newline

\par We then claim that knowing
$\Sigma_{u\in U^{0,0,A_1,0,\frac{n}{2}-A_1}} a_u
\\ Tail^{Shad}_{+} [C^u_{g^n}(\phi)]$ we can determine
$\Sigma_{u\in U^{0,0,A_1,0,\frac{n}{2}-A_1}} a_u
 C^u_{g^n}(\phi)$. Specifically, we show that for the complete contraction
$C^{+}_{g^n}(\phi)=C^u_{g^n}(\phi)$, $u\in
U^{0,0,A_1,0,\frac{n}{2}-A_1}$, we have that:

\begin{equation}
\label{arv}
Tail^{Shad}_{+}[C^{+}_{g^n}(\phi)]= (-1)^{\frac{n}{2}-A_1-1}2A_1
\cdot (\frac{n}{2}-A_1)\cdot C^\sharp_{g^n}(\phi,\vec{\xi})
\end{equation}

{\it Proof of (\ref{arv}):}
\newline

\par For any $\vec{\xi}$-contraction $C_{g^n}(\phi,\vec{\xi})$, we
 denote by $O^{Shad}_{+}[C_{g^n}(\phi,\vec{\xi})]$ the sublinear
 combination in $O^{Shad}[C_{g^n}(\phi,\vec{\xi})]$ of
$\vec{\xi}$-contractions in the form (\ref{stewart}).

 Firstly, we denote by $C^{+}_{g^n}(\phi,
\vec{\xi})$ the descendant of $C^{+}_{g^n}(\phi)$ that arises by
 replacing each of the $A_1$ factors $R$ by $|\vec{\xi}|^2$ and each
of the $\frac{n}{2}-A_1$ factors $\Delta\phi$ by
$\vec{\xi}^i\nabla_i\phi$ (in the $N$-cancelled notation).
 We observe that for any descendant $C'_{g^n}(\phi,\vec{\xi})$
of $C^{+}_{g^n}(\phi)$ other than the above, we will have that
$Tail^{Shad}_{+}[C'_{g^n}(\phi,\vec{\xi})]$=0. This is true by virtue
 of the same arguments as for the previous case (at $\vec{\xi}$-length
$\frac{n}{2}$ there must be an internal contraction). Hence, it suffices
 to show that $O^{Shad}_{+}[C^{+}_{g^n}(\phi,\vec{\xi})]$
is equal to the right hand side of (\ref{arv}).

\par So, let us begin by performing the iterative integration by parts.
 We first integrate by parts the factor $\vec{\xi}^i$ that
contracts against the first factor $\nabla_i\phi$. Note that, although
we have imposed restrictions on the order of our integrations by parts,
 in this case we can pick an order so that we first integrate by parts
with respect to this factor $\vec{\xi}$. If $\nabla_i$ hits a factor
$\nabla\phi$ or a factor $\vec{\xi}$ that does not contract against
 another factor $\vec{\xi}$, we denote the
$\vec{\xi}$-contraction that is generically thus obtained  by
$C^d_{g^n}(\phi,\vec{\xi})$. We observe that
$O^{Shad}_{+}[C^d_{g^n}(\phi,\vec{\xi})]=0$, since each
$\vec{\xi}$-contraction in that sublinear combination will either have
 length $\ge\frac{n}{2} +1$ or at least one factor $\nabla^a\phi$,
$a\ge 2$.
If $\nabla_i$ hits a factor $|\vec{\xi}|^2$, we obtain an expression
$2\nabla_i\vec{\xi}_j\vec{\xi}^j$ and we denote the
$\vec{\xi}$-contraction that we have obtained by
$C^{*}_{g^n}(\phi,\vec{\xi})$.
We then proceed to integrate by parts the factor $\vec{\xi}^j$.

\par Now, if $\nabla^j$ hits a factor $\vec{\xi}$ that does not
contract against another factor $\vec{\xi}$ or if it hits a factor
$|\vec{\xi}|^2$, we
generically denote the $\vec{\xi}$-contraction that is thus obtained by
$C^d_{g^n}(\phi,\vec{\xi})$ and we observe that $O^{Shad}_{+}
[C^d_{g^n}(\phi,\vec{\xi})]=0$. This follows because in the first
case we will obtain a $\vec{\xi}$-contraction of
$\vec{\xi}$-length
$\ge\frac{n}{2}+1$ and in the second we will have less than $A_1-1$
 factors $|\vec{\xi}|^2$.

\par Initially, we consider the $\vec{\xi}$-contraction $C^{*,1}_{g^n}(\phi,\vec{\xi})$
that arises when  $\nabla^j$ hits the first factor $\nabla_i\phi$.
In that case, $C^{*,1}_{g^n}(\phi,\vec{\xi})$ is the
complete contraction:

\begin{equation}
\label{stewart2}
contr((|\vec{\xi}|^2)^{A_1-1}\otimes\nabla^2_i\vec{\xi}_j\otimes
\nabla^{ij}\phi\otimes
(\nabla)^{r_1}\phi\vec{\xi}_{r_1}\otimes\dots \otimes
(\nabla)^{r_{\frac{n}{2}-A_1-1}}\phi\vec{\xi}_{r_{\frac{n}{2}-A_1-1}})
\end{equation}

We show that $O^{Shad}_{+}[C^1_{g^n}(\phi,\vec{\xi})]=
(-1)^{\frac{n}{2}-A_1-1}C^\sharp_{g^n}(\phi,\vec{\xi})$.
\newline

\par This follows by the iterative integrations by parts procedure.
 The algorithm to obtain $(-1)^{\frac{n}{2}-A_1-1}C^\sharp_{g^n}(\phi,\vec{\xi})$
is to successively integrate by parts each of the
$\frac{n}{2}-A_1-1$ factors $\vec{\xi}$ that contract against a
factor $\nabla\phi$ and make it hit the one factor
$S\nabla^p\vec{\xi}$ and then symmetrize. We then obtain
$(-1)^{\frac{n}{2}-A_1-1}C^\sharp_{g^n}(\phi,\vec{\xi})$. We
observe that if at any stage we integrate by parts a factor
$\vec{\xi}$ and hit the factor $\nabla^2\phi$ or a factor
$|\vec{\xi}|^2$ or a factor $\vec{\xi}$ or a factor $\nabla\phi$,
then performing the rest of the iterative integrations by parts we
will not obtain a $\vec{\xi}$-contraction in the form
$C^\sharp_{g^n}(\phi,\vec{\xi})$.
\newline

\par On the other hand, we consider the $\vec{\xi}$-contraction that
 arises when $\nabla^j$ hits the $h^{th}$ factor $\nabla\phi$,
$h\ge 2$. We denote the $\vec{\xi}$-contraction that arises thus
by $C^{*,h}_{g^n}(\phi,\vec{\xi})$. We then claim that
$O^{Shad}_{+}[C^{*,h}_{g^n}(\phi,\vec{\xi})]=
(-1)^{\frac{n}{2}-A_1-1}C^\sharp_{g^n}(\phi,\vec{\xi})$. It is
clear that if we can show the
 above claim, (\ref{arv}) will follow immediately.
\newline

\par To see this, we initially observe that up to permuting
factors $\nabla\phi$, $C^{*,h}_{g^n}(\phi,\vec{\xi})$ is in the
form:

\begin{equation}
\label{stewart50}
contr((|\vec{\xi}|^2)^{A_1-1}\otimes\nabla_i\vec{\xi}_j\otimes
(\nabla)^i\phi\otimes
(\nabla)^j(\nabla)^{r_1}\phi\vec{\xi}_{r_1}\otimes\dots \otimes
(\nabla)^{r_{\frac{n}{2}-A_1-1}}\phi\vec{\xi}_{r_{\frac{n}{2}-A_1-1}})
\end{equation}

\par Moreover, it follows that
$(-1)^{\frac{n}{2}-A_1-1}C^\sharp_{g^n}(\phi,\vec{\xi})$ arises in
$O^{Shad}_{+}[C^{*,h}_{g^n}(\phi,\vec{\xi})]$ when we integrate by
parts all the factors $\vec{\xi}_a$ and hit the factor
$S\nabla\vec{\xi}$ and then replace
$\nabla^{\frac{n}{2}-A_1-1}\nabla\vec{\xi}$ by
$S\nabla^{\frac{n}{2}-A_1}\vec{\xi}$. We observe that if we
perform any other integration by parts, we will not obtain
$C^\sharp_{g^n}(\phi,\vec{\xi})$: If we hit a factor
$|\vec{\xi}|^2$ by a $\nabla$, we will obtain a
$\vec{\xi}$-contraction with fewer than $A_1-1$ factors
$|\vec{\xi}|^2$. If we hit a factor $\vec{\xi}$ that does not
contract against another factor $\vec{\xi}$, we will have
$\vec{\xi}$-length $\ge \frac{n}{2}+1$. If we hit a factor
$\nabla\phi$ or the factor $\nabla^2\phi$, we will respectively
have two factors $S\nabla^p\phi$ with $p\ge 2$ or one factor
$S\nabla^p\phi$ with $p\ge 3$. Finally, if we hit the factor
$S\nabla^p\vec{\xi}$ by a derivative $\nabla_i$ and
anti-symmetrize using the equation:

\begin{equation}
\label{symunsymgenxi}
\begin{split}
&{\nabla}_aS{\nabla}^m_{r_1\dots r_m}\vec{\xi}_j=
S{\nabla}^m_{ar_1\dots r_m}\vec{\xi}_j +
C_{m-1}\cdot
S^{*}{\nabla}^{m-1}_{r_1\dots r_{m-1}}
R_{aijd}\vec{\xi}^d +
\\&
{\Sigma}_{u\in U^m} a_u
pcontr({\nabla}^{m'} R_{abcd} S{\nabla}^{s_u}\vec{\xi})
\end{split}
\end{equation}
from \cite{a:dgciI} (and the notational conventions there),
we obtain a $\vec{\xi}$-contraction with
 a factor of the form $\nabla^mR_{ijkl}$. Hence, by the iterative
 integrations by parts procedure, the $O^{Shad}$ of such a factor will
 consist of $\vec{\xi}$-contractions with a factor $\nabla^mR_{ijkl}$,
 so we have completely shown our claim.

\par Hence, we have determined
$\Sigma_{u\in U^{0,0,A_1,0,\frac{n}{2}-A_1}} a_u Tail^{Shad}_{+}
[C^u_{g^n}(\phi)]$. In other words, we have determined the constant
$(Const)_\sharp$ for which:

\begin{equation}
\label{xysia}
\Sigma_{u\in U^{0,0,A_1,0,\frac{n}{2}-A_1}} a_u Tail^{Shad}_{+}
[C^u_{g^n}(\phi)]=(Const_\sharp)\cdot C^\sharp_{g^n}(\phi,\vec{\xi})
\end{equation}
Now, we only have to replace each expression
$|\vec{\xi}|^2$ by an expression $R$ and the expression

$$\nabla^{\frac{n}{2}-A_1+1}_{r_1\dots
r_{\frac{n}{2}-A_1-1}i}\vec{\xi}_j\otimes \nabla^{ij}\phi\otimes
(\nabla)^{r_1}\phi\otimes\dots \otimes
(\nabla)^{r_{\frac{n}{2}-A_1-1}}\phi$$ by an expression $R\cdot
(\Delta\phi)^{\frac{n}{2}-A_1}$. Then, using
 (\ref{arv}) and (\ref{xysia}), we determine the constant $Const'$ for
which:

$$\Sigma_{u\in U^{0,0,A_1,0,\frac{n}{2}-A_1}} a_u C^u_{g^n}(\phi)=
(Const')\cdot R^{A_1}\cdot (\Delta\phi)^{\frac{n}{2}-A_1}$$

In other words, we determine the sublinear combination
$\Sigma_{u\in U^{0,0,A_1,0,\frac{n}{2}-A_1}} a_u C^u_{g^n}(\phi)$.
That concludes the proof of our second claim.
\newline

\subsection{Determining the sublinear combination
\\$\Sigma_{u\in
U^{A_1-X_1-C_1,X_1,C_1,\frac{n}{2}-\Delta_1,\Delta_1}} a_u
C^u_{g^n}(\phi)$.}

\par We call the list $(A_1-X_1-C_1,X_1,C_1,\frac{n}{2}-A_1
-\Delta_1,\Delta_1)$ the {\it critical list}. We denote the
 index set
$U^{A_1-X_1-C_1,X_1,C_1,\frac{n}{2}-A_1-\Delta_1,\Delta_1}$ by
$U^{crit}$ for short. Moreover, whenever we refer to a list
$(Z,X,C,\Gamma,\Delta)$ for which we have not yet determined
$\Sigma_{u\in U^{Z,X,C,\Gamma,\Delta}} a_u C^u_{g^n}(\phi)$, we
will say that the list $(Z,X,C,\Gamma,\Delta)$ is {\it subsequent}
 to the critical list. We will also say that $u$ or
$C^u_{g^n}(\phi)$ is subsequent to the critical list when
$u\in U^{Z,X,C,\Gamma,\Delta}$.

On the other hand, for each list
$(Z,X,C,\Gamma,\Delta)$ where we have determined $\Sigma_{u\in
U^{Z,X,C,\Gamma,\Delta}} a_u C^u_{g^n}(\phi)$, we will say that
the list $(Z,X,C,\Gamma,\Delta)$ {\it preceded} the critical
 list. Accordingly, in that case,
if $u\in U^{Z,X,C,\Gamma,\Delta}$, we will say that
$u$ or $C^u_{g^n}(\phi)$ preceded the critical list.

\par  We will distinguish three cases and separately
 prove our claim in each of those cases. The first case is when
$\Delta_1 <\frac{n}{2}-A_1$. The second one is when
$\Delta_1=\frac{n}{2}-A_1$ and $X_1>0$. The third is when
$\Delta_1=\frac{n}{2}-A_1$, $X_1=0$. In the third case we observe that
 we will have that $X_1+C_1<A_1$ (otherwise we are in the base case
 that we have already dealt with). In each of the three cases,
we will use the equation:

\begin{equation}
\label{manfr3}
I^{\frac{n}{2}-A_1}_{g^n}(\phi)= \Sigma_{k\in K} a_k C^k_{g^n}(\phi)+
\Sigma_{u\in U} a_u C^u_{g^n}(\phi)
\end{equation}
which holds modulo complete contractions of length
$\ge\frac{n}{2}+1$. We recall that the sublinear combination
$\Sigma_{k\in K} a_k C^k_{g^n}(\phi)$ is known, and each sublinear
combination $U^{Z,X,C,\Gamma,\Delta}$, where $(Z,X,C,\Gamma,\Delta)$
 precedes $U^{crit}$ is also known.

\par We proceed to prove our claim in each of the three cases.
\newline

{\it The first case.} We consider
$Shad[I^{\frac{n}{2}-A_1}_{g^n}(\phi)]$ and focus on the sublinear
 combination of $\vec{\xi}$-contractions in the following form:

\begin{equation}
\label{arv2}
\begin{split}
&contr(R_{i_1j_1k_1l_1}\otimes\dots\otimes R_{i_{A_1-X_1-C_1}
j_{A_1-X_1-C_1}k_{A_1-X_1-C_1}l_{A_1-X_1-C_1}}\otimes
\nabla_{a_1}\vec{\xi}_{b_1}\otimes\dots\otimes
\\&\nabla_{a_{X_1}}
\vec{\xi}_{b_{X_1}}
\otimes(|\vec{\xi}|^2)^{C_1}\otimes
S(\nabla^{\Delta_1}_{s_1\dots s_{\Delta_1}}\nabla^2_{f_1g_1})\phi
\otimes\nabla^2_{f_2g_2}\phi \otimes \dots\otimes
\nabla^2_{f_{\frac{n}{2}-A_1-\Delta_1}
g_{\frac{n}{2}-A_1-\Delta_1}}\phi
\\&\otimes (\nabla)^{s_1}\phi\otimes\dots \otimes
(\nabla)^{s_{\Delta_1}}\phi)
\end{split}
\end{equation}

\par We denote the sublinear combination of
$\vec{\xi}$-contractions in the form (\ref{arv2}) in
$Shad[I^{\frac{n}{2}-A_1}_{g^n}(\phi)]$ by
$Shad_o[I^{\frac{n}{2}-A_1}_{g^n}(\phi)]$. We then claim that:

\begin{equation}
\label{sisy}
Shad_o[I^{\frac{n}{2}-A_1}_{g^n}(\phi)]=0
\end{equation}

\par This can be seen by the following reasoning: We write out
the sublinear combination of $\vec{\xi}$-contractions of
$\vec{\xi}$-length $\frac{n}{2}$ in
$Shad[I^{\frac{n}{2}-A_1}_{g^n}(\phi)]$ as a linear combination of
$\vec{\xi}$-contractions in the form:

\begin{equation}
 \label{stigma}
\begin{split}
&contr({\nabla}_{r_1\dots r_{m_1}}^{m_1}R_{i_1j_1k_1l_1}\otimes
\dots \otimes {\nabla}_{v_1\dots
v_{m_s}}^{m_s}R_{i_sj_sk_sl_s}\otimes \\& {\nabla}_{t_1 \dots
t_{p_1}}^{p_1}Ric_{{\alpha}_1 {\beta}_1} \otimes \dots \otimes
{\nabla}_{z_1 \dots z_{p_q}}^{p_q} Ric_{{\alpha}_q{\beta}_q}
\otimes{\nabla}^{{\nu}_1}_ {{\chi}_1\dots
{\chi}_{{\nu}_1}}\phi \otimes\dots \otimes
{\nabla}^{{\nu}_Z}_{{\omega}_1\dots {\omega}_{{\nu}_Z}}
\phi
\\& \otimes
S{\nabla}^{{\mu}_1}\vec{\xi}_{j_1}\dots \dots
S{\nabla}^{{\mu}_r}\vec{\xi}_{j_s}\otimes |\vec{\xi}|^2
\otimes \dots \otimes |\vec{\xi}|^2)
\end{split}
\end{equation}

with $Z=\frac{n}{2}-A_1$. Then, we define $Tail^{Shad}_\alpha
[I^{\frac{n}{2}-A_1}_{g^n}(\phi)]$ to stand for
 the sublinear combination in
$Tail^{Shad} [I^{\frac{n}{2}-A_1}_{g^n}(\phi)]$ that consists of
$\vec{\xi}$-contractions of $\vec{\xi}$-length $\frac{n}{2}$ for
which the decreasing rearrangement of the list $\nu_1,\dots
,\nu_{\frac{n}{2}-A_1}$ is $(\Delta_1+2,2,\dots , 2,1,\dots ,1)$
(we are writing the number $2$ $\Gamma_1 -1$ times and $1$
$\Delta_1$ times. Then, by Lemma \ref{firstsymcanc}, we have that:

\begin{equation}
\label{sisy2}
Shad_\alpha[I^{\frac{n}{2}-A_1}_{g^n}(\phi)]=0
\end{equation}

\par Now, we consider the sublinear combination
$Shad_{\alpha,\beta}[I^{\frac{n}{2}-A_1}_{g^n}(\phi)]$ in
\\ $Shad_\alpha
[I^{\frac{n}{2}-A_1}_{g^n}(\phi)]$ where there
are no factors with internal contractions (in particular there are
no factors $\nabla^pRic$ or $\nabla^mR_{ijkl}$ with internal
contractions). Then, since the number of internal contractions
remains invariant under the permutations of definition
7 in \cite{a:dgciI}, modulo introducing
$\vec{\xi}$-contractions of
$\vec{\xi}$-length $\ge\frac{n}{2}+1$, we will have that modulo
$\vec{\xi}$-contractions of $\vec{\xi}$-length $\ge\frac{n}{2}+1$:

\begin{equation}
\label{sisy2}
Shad_{\alpha,\beta}[I^{\frac{n}{2}-A_1}_{g^n}(\phi)]=0
\end{equation}

\par Moreover, we define $Shad_{\alpha,\beta,\gamma}
[I^{\frac{n}{2}-A_1}_{g^n}(\phi)]$ to stand for the sublinear
 combination in $Shad_{\alpha,\beta}[I^{\frac{n}{2}-A_1}_{g^n}(\phi)]$
where the $\Delta_1$ factors $\nabla\phi$ are all contracting against
 the one factor $\nabla^{\Delta_1+2}\phi$. We observe that
the number of factors $\nabla\phi$ that contract against the
factor $\nabla^{\Delta_1+2}\phi$ remains invariant under the
permutations allowed by definition 7 in \cite{a:dgciI}, modulo
introducing $\vec{\xi}$-contractions of $\vec{\xi}$-length
$\ge\frac{n}{2}+1$. Hence, we have that modulo
$\vec{\xi}$-contractions of $\vec{\xi}$-length $\ge\frac{n}{2}+1$:

\begin{equation}
\label{sisy3}
Shad_{\alpha,\beta,\gamma}[I^{\frac{n}{2}-A_1}_{g^n}(\phi)]=0
\end{equation}

\par Finally, we define
$Shad_{\alpha,\beta,\gamma,\delta}[I^{\frac{n}{2}-A_1}_{g^n}(\phi)]$
to stand for the sublinear combination in
$Shad_{\alpha,\beta,\gamma}[I^{\frac{n}{2}-A_1}_{g^n}(\phi)]$ that
consists of the $\vec{\xi}$-contractions with $X_1$ factors
$\nabla\vec{\xi}$ and no more factors of the form
$S\nabla^u\vec{\xi}$ and, in addition, with $C_1$ factors
$|\vec{\xi}|^2$. Since both the number of factors
$S\nabla^p\vec{\xi}$ ($p\ge 1$) and the number of such factors for
which $p=1$, and also the number of factors $|\vec{\xi}|^2$ is
invariant under the permutations of
 definition 7 in \cite{a:dgciI}, we have
 that modulo $\vec{\xi}$-contractions of $\vec{\xi}$-length
$\ge\frac{n}{2}+1$:

\begin{equation}
\label{sisy3}
Shad_{\alpha,\beta,\gamma,\delta}[I^{\frac{n}{2}-A_1}_{g^n}(\phi)]=0
\end{equation}

\par Now, we observe that
$Shad_{\alpha,\beta,\gamma,\delta}[I^{\frac{n}{2}-A_1}_{g^n}(\phi)]$
indeed consists of $\vec{\xi}$-contractions of the form
(\ref{arv2}). This follows just because we are considering
$\vec{\xi}$-length $\frac{n}{2}$ and weight $-n$. Hence, we must
have $A_1-X_1-C_1$ factors $\nabla^mR_{ijkl}$ with no internal
contractions. But since the $\vec{\xi}$-contractions in the form
(\ref{arv2}) have indeed weight $-n$, it follows that any
$\vec{\xi}$-contraction with the restrictions above and with at
least one factor $\nabla^mR_{ijkl}, m>0$ cannot have weight $-n$.
\newline

\par Now, for each complete contraction $C_{g^n}(\phi)$
in $I^{\frac{n}{2}-A_1}_{g^n}(\phi)$, we denote by
$Tail^{Shad}_o [C_{g^n}(\phi)]$ the sublinear combination
of $\vec{\xi}$-contractions in the form  (\ref{arv2}) in
$Tail^{Shad}[C_{g^n}(\phi)]$. This notation extends to linear
 combinations.

\par Now, if we write  $I^{\frac{n}{2}-A_1}_{g^n}(\phi)$ out as in
 (\ref{manfr3}), we claim that for each $C^u_{g^n}(\phi)$, where $u$ is
 subsequent to the critical list, we have that modulo
$\vec{\xi}$-contractions of $\vec{\xi}$-length $\ge\frac{n}{2}+1$:

\begin{equation}
\label{althu}
Tail^{Shad}_o [C^u_{g^n}(\phi)]=0
\end{equation}

\par We will prove this below. For now, we note how
we can then determine our desired sublinear combination
$\Sigma_{u\in U^{crit}} a_u C^u_{g^n}(\phi)$.
 Initially we observe that
if we can show (\ref{althu}), we will then be able to determine the
 sublinear combination $\Sigma_{u\in U^{crit}} a_u Tail^{Shad}_o[C^u_{g^n}(\phi)]$
from equation (\ref{sisy}). We then also claim that for each
$u\in U^{crit}$, the sublinear combination
$Tail^{Shad}_o [C^u_{g^n}(\phi)]$ is obtained from $C^u_{g^n}(\phi)$
by performing the following algorithm:
We replace each
 factor $R$ by $-|\vec{\xi}|^2$, each factor $Ric_{ij}$ by
$-\nabla_i\vec{\xi}_j$ and each factor $\Delta\phi$ by
$\vec{\xi}^i\nabla_i\phi$ (in $N$-cancelled notation). We then
 integrate by parts the $\Delta_1$ factors $\vec{\xi}$ that contract
 against factors $\nabla\phi$ and make each $\nabla_i$ that arises thus
hit the same factor $\nabla^2\phi$.

\par This follows just by the iterative
 integration by parts procedure, and the same arguments as above.
  Since we have determined
$\Sigma_{u\in U^{crit}} a_u Tail^{Shad}_o[C^u_{g^n}(\phi)]$, then
by replacing each expression $|\vec{\xi}|^2$ by $R$, each
expression $\nabla_i\vec{\xi}$ by $-Ric_{ij}$ and each expression
$\nabla^{\Delta_1}_{s_1\dots
s_{\Delta_1}}(\nabla^2_{f_1g_1})\phi\otimes \dots\otimes
\nabla^2_{f_{\frac{n}{2}-A_1-\Delta_1}
g_{\frac{n}{2}-A_1-\Delta_1}}\phi\otimes
(\nabla)^{s_1}\phi\otimes\dots \otimes
(\nabla)^{s_{\Delta_1}}\phi)$ by $(\nabla^2_{f_1g_1}\phi\otimes
\dots\otimes \nabla_{f_{\frac{n}{2}-A_1-\Delta_1}
g_{\frac{n}{2}-A_1-\Delta_1}}\phi)(\Delta\phi)^{\Delta_1}$, we
 have determined the sublinear combination $\Sigma_{u\in
U^{crit}} a_u C^u_{g^n}(\phi)$. Moreover, we see that by construction,
 the pattern of those particular contractions between indices in
factors $R_{ijkl}, Ric_{ij}, \nabla^2\phi$ is preserved.
\newline

\par So, matters are reduced to showing that for each
$C^u_{g^n}(\phi)$  where $u$ is subsequent to the critical character,
we must have that $Tail^{Shad}_o [C^u_{g^n}(\phi)]=0$.
Firstly, we observe that we may restrict attention to the
descendants of $C^u_{g^n}(\phi)$ that do not have internal
contractions. This follows by the same reasoning as in the previous
 case. Then, we observe that if $C^u_{g^n}(\phi)$ has
$\Delta<\Delta_1$ factors $\Delta\phi$, then each
$\vec{\xi}$-contraction
 of length $\frac{n}{2}$ in $Tail^{Shad}[C^u_{g^n}(\phi)]$
will have less than $\Delta_1$ factors $\nabla\phi$.  Similarly,
if $C^u_{g^n}(\phi)$ has less than $C_1$ factors $R$ then each
 complete contraction of $\vec{\xi}$-length
$\frac{n}{2}$ in $Tail^{Shad}[C^u_{g^n}(\phi)]$ will have
 less than $C_1$ expressions $|\vec{\xi}|^2$.
Finally, if $C^u_{g^n}(\phi)$
has $\Delta_1$ factors $\Delta\phi$, $C_1$ factors $R$ and
 less than $X_1$ factors $Ric_{ij}$, then each $\vec{\xi}$-contraction
 in $Tail^{Shad}[C^u_{g^n}(\phi)]$
will either have less than $X_1$ factors $\nabla_i\vec{\xi}$ or
 less than $C_1$ expressions $|\vec{\xi}|^2$. Thus we have shown our
 claim.
\newline

\par {\it The second case}, where $\Delta_1=\frac{n}{2}-A_1$
and $X_1>0$. We again consider the shadow divergence formula for
$I^{\frac{n}{2}-A_1}_{g^n}(\phi)$, and we focus on the sublinear
 combination of $\vec{\xi}$-contractions in the form:

\begin{equation}
\label{arv4}
\begin{split}
&contr(R_{i_1j_1k_1l_1}\otimes\dots\otimes R_{i_{A_1-X_1-C_1}
j_{A_1-X_1-C_1}k_{A_1-X_1-C_1}l_{A_1-X_1-C_1}}\otimes
\\& S\nabla^{\frac{n}{2}-A_1+1}_{s_1\dots s_{\frac{n}{2}-A_1}a_1}
\vec{\xi}_{b_1}\otimes\nabla_{a_2}\vec{\xi}_{b_2}\otimes
\dots\otimes  \nabla_{a_{X_1}}
\vec{\xi}_{b_{X_1}}\otimes(|\vec{\xi}|^2)^{C_1} \otimes
(\nabla)^{s_1}\phi\otimes\dots \otimes
(\nabla)^{s_{\Delta_1}}\phi)
\end{split}
\end{equation}

We denote the above sublinear combination by
$Shad_{+}[I^{\frac{n}{2}-A_1}_{g^n}(\phi)]$. Since the
 shadow divergence formula holds formally, by an analogous argument as
for the previous case, it follows that:

\begin{equation}
\label{arv5}
Shad_{+}[I^{\frac{n}{2}-A_1}_{g^n}(\phi)]=0
\end{equation}
For each complete contraction $C_{g^n}(\phi)$ in
$I^{\frac{n}{2}-A_1}_{g^n}(\phi)$, we denote by $Tail^{Shad}_o
[C^u_{g^n}(\phi)]$ the sublinear combination of
$\vec{\xi}$-contractions in the form (\ref{arv4}) in
$Tail^{Shad}[C_{g^n}(\phi)]$. (This is not the same as the
previous $Tail^{Shad}[C_{g^n}(\phi)]$).
\newline

\par Now, by a similar reasoning as for the previous case, we observe
that for each $C^u_{g^n}(\phi)$ that is subsequent to the critical
 character we have
$Tail^{Shad}_o[C^u_{g^n}(\phi)]=0$. This follows because if
$C^u_{g^n}(\phi)$ has either less than $\Delta_1$ factors
$\Delta\phi$, or $\Delta_1$ such factors and less than $C_1$
factors $R$ or $C_1$ such factors and less than $X_1$ factors
$Ric$. In those cases, we respectively have that each
$\vec{\xi}$-contraction in $Tail[C^u_{g^n}(\phi)]$ will have less
than $\Delta_1$ factors $\nabla\phi$ or less than $C_1$ factors
$|\vec{\xi}|^2$ or less than $X_1$ factors
$S\nabla^p\vec{\xi}$. Hence, using (\ref{arv5}), we determine the sublinear
combination $\Sigma_{u\in U^{crit}} a_u
Tail^{Shad}_o[C^u_{g^n}(\phi)]$.
\newline

\par We now claim that for each $C^u_{g^n}(\phi), u\in U^{crit}$, the
sublinear combination $Tail^{Shad}_o[C^u_{g^n}(\phi)]$ arises as
follows: We initially replace each of the $C_1$ factors $R$ by
$|\vec{\xi}|^2$, each of the $X_1$ factors $Ric_{ij}$ by
$-\nabla_i\vec{\xi}_j$ and each of the $\frac{n}{2}-A_1$ factors
$\Delta\phi$ by $\nabla^i\phi\vec{\xi}_i$ (we are using
$N$-cancelled notation). We then integrate by parts the
$\frac{n}{2}-A_1$ factors $\vec{\xi}$ that contract against a
factor $\nabla\phi$ and make the derivatives $\nabla^i$ hit the
same one factor $\nabla_i\vec{\xi}_j$ and replace
$\nabla^{\frac{n}{2}-A_1}_{i_1\dots i_{\frac{n}{2}-A_1}}
\nabla_i\vec{\xi}_j$ by $S\nabla^{\frac{n}{2}-A_1}_{i_1\dots
i_{\frac{n}{2}-A_1}i}\vec{\xi}_j$. This follows by the iterative
integrations by parts procedure, as in the previous case.

\par Therefore, once we have determined
$\Sigma_{u\in U^{crit}} a_u Tail^{Shad}_o[C^u_{g^n}(\phi)]$, we
can determine $\Sigma_{u\in U^{crit}} a_u C^u_{g^n}(\phi)$ as
follows: We replace each factor $|\vec{\xi}|^2$ by $R$, each
factor $\nabla_i\vec{\xi}_j$ by $-Ric_{ij}$ and each expression
$S\nabla^{\frac{n}{2}-A_1+1}_{s_1\dots s_{\frac{n}{2}-A_1}a_1}
\vec{\xi}_{b_1}(\nabla)^{s_1}\phi\otimes\dots \otimes
(\nabla)^{s_{\Delta_1}}\phi$ by $Ric_{a_1b_1}
(\Delta\phi)^{\frac{n}{2}-A_1}$. We then determine the sublinear
combination $\Sigma_{u\in U^{crit}} a_u C^u_{g^n}(\phi)$.
\newline

{\it The third case.}
\newline

\par Finally, we have to consider the third case. We now consider
$I^{\frac{n}{2}-A_1}_{g^n}(\phi)$ and distinguish the two subcases
$C_1=0$ or $C_1>0$.
\newline

{\it The first subcase $C_1=0$}. Modulo complete contractions of
length $\ge\frac{n}{2}+1$, we write out
$I^{\frac{n}{2}}_{g^n}(\phi)$ in the form:

\begin{equation}
\label{ishbilia} I^{\frac{n}{2}}_{g^n}(\phi)= \Sigma_{g\in G} a_g
C^g_{g^n}(\phi)+\Sigma_{u\in U^{crit}} a_u
C^u_{g^n}(\phi)+\Sigma_{u\in U^{subs}}a_u C^u_{g^n}(\phi)
\end{equation}
where $\Sigma_{g\in G} a_g C^g_{g^n}(\phi)$ stands for the known
sublinear combination in $I^{\frac{n}{2}}_{g^n}(\phi)$ (this now
includes a part of $\Sigma_{u\in U} a_u C^u_{g^n}(\phi)$).
$\Sigma_{u\in U^{crit}} a_u C^u_{g^n}(\phi)$ stands for the
sublinear combination of compete contractions indexed in the
critical list, $U^{crit}$. Finally, $\Sigma_{u\in U^{subs}}a_u
C^u_{g^n}(\phi)$ stands for the sublinear combination of complete
contractions $C^u_{g^n}(\phi)$ that are subsequent to the critical
list.

\par We focus on
the super divergence formula for
$I^{\frac{n}{2}}_{g^n}(\phi)$. We
pick out the sublinear combination of  complete contractions  in
the form:

\begin{equation}
\label{arv6} contr(\nabla^{\Delta_1}_{r_1\dots r_{\Delta_1}}
R_{i_1j_1k_1l_1}\otimes R_{i_2j_2k_2l_2}\otimes \dots \otimes
 R_{i_{A_1}j_{A_1}k_{A_1}l_{A_1}}
\otimes(\nabla)^{s_1}\phi\otimes\dots\otimes
(\nabla)^{s_{\frac{n}{2}-A_1}} \phi)
\end{equation}
where each of the factors $\nabla\phi$ contracts against an index
in the factor $\nabla^{\Delta_1}R_{ijkl}$.

\par We denote the corresponding sublinear combination of complete
contractions in $supdiv[I^{\frac{n}{2}-A_1}_{g^n}(\phi)]$
by $supdiv_{+}[I^{\frac{n}{2}-A_1}_{g^n}(\phi)]$. Since the super
divergence formula holds formally, it follows that:

$$supdiv_{+}[I^{\frac{n}{2}-A_1}_{g^n}(\phi)]=0$$
modulo complete contractions of length $\ge\frac{n}{2} +1$. Now,
for each $C_{g^n}(\phi)$ in $I^{\frac{n}{2}-A_1}_{g^n}(\phi)$, we
denote by $Tail_{+}[C_{g^n}(\phi)]$ the sublinear combination in
$Tail[C_{g^n}(\phi)]$ that consists of complete contractions in
the form (\ref{arv6}).

\par We then again observe that for each $u$ that is subsequent
to the critical list, we have $Tail_{+}[C^u_{g^n}(\phi)]=0$. This
follows since if $C^u_{g^n}(\phi)$ is subsequent to the critical
list it must have less than $\frac{n}{2}-A_1$ factors
$\Delta\phi$, hence any complete contraction of length
$\frac{n}{2}$ in $Tail[C^u_{g^n}(\phi)$ must have less than
$\frac{n}{2}-A_1$ factors $\nabla\phi$. On the other hand, for
each $u\in U^{crit}$ we have that $Tail_{+}[C^u_{g^n}(\phi)]$
arises from $C^u_{g^n}(\phi)$ as
 follows: We replace each of the factors $\Delta\phi$ by
$\nabla_i\phi\vec{\xi}^i$ and then integrate by parts the
$\frac{n}{2}-A_1$ factors $\vec{\xi}$ and make each of them hit the
same factor $R_{ijkl}$ (there are $A_1$ choices of the factor
$R_{ijkl}$ that we may pick).
 The sublinear combination that arises thus is $Tail_{+}[C^u_{g^n}(\phi)]$.
In fact, we observe that if $C^u_{g^n}(\phi)$ is of the form:

\begin{equation}
\label{paisley} contr(R_{i_1j_1k_1l_1}\otimes\dots \otimes
R_{i_{A_1}j_{A_1}k_{A_1}l_{A_1}}\otimes
(\Delta\phi)^{\frac{n}{2}-A_1})
\end{equation}

Then $Tail_{+}[C^u_{g^n}(\phi)]$ can be written as a sum of $A_1$
complete contractions in the form:

\begin{equation}
\label{kati} (-1)^{\frac{n}{2}-A_1}
contr(\nabla^{\frac{n}{2}-A_1}_{i_1\dots i_{\frac{n}{2}-A_1}}
R_{ijkl}\otimes\dots \otimes R_{i'j'k'l'} \otimes
(\nabla)^{i_1}\phi\otimes\dots\otimes
(\nabla)^{i_{\frac{n}{2}-A_1}} \phi)
\end{equation}
where the $h^{th}$ term in the sum arises from $C^u_{g^n}(\phi)$
by replacing all the factors $\Delta\phi$ by a factor
$\nabla_{a_j}\phi$ ($1\le j\le \frac{n}{2}-A_1$) and then hitting
the $h^{th}$ factor $R_{ijkl}$ in $C^u_{g^n}(\phi)$ by
$\frac{n}{2}-A_1$ derivatives $(\nabla)^{a_j}$. In order to
facilitate our work further down, we will write out:

\begin{equation}
\label{clover} Tail_{+}[C^u_{g^n}(\phi)]=\Sigma_{h=1}^{A_1}
C^{u,h}_{g^n}(\phi)
\end{equation}
where $C^{u,h}_{g^n}(\phi)$ stands for the $h^{th}$ complete
contraction explained above. Given the form (\ref{paisley}) of
$C^u_{g^n}(\phi)$, we have that $C^{u,h}_{g^n}(\phi)$ will be in
the form:

\begin{equation}
\label{arv12} contr(R_{i_1j_1k_1l_1}\otimes\dots \otimes
\nabla_{i_1\dots i_{\Delta_1}}R_{i_hj_hk_hl_h}\otimes\dots \otimes
R_{i_{A_1}j_{A_1}k_{A_1}l_{A_1}}\otimes
(\nabla)^{s_1}\phi\otimes\dots \otimes
(\nabla)^{s_{\Delta_1}}\phi)
\end{equation}

\par Now, for each $u\in U^{crit}$, we denote by
$C^u(g^n)$ the complete contraction of weight $-2A_1$:

$$contr(R_{i_1j_1k_1l_1}\otimes\dots \otimes R_{i_{A_1}j_{A_1}k_{A_1}
l_{A_1}})$$

\par We then claim that we can determine the linear combination
$\Sigma_{u\in U^{crit}} a_u C^u(g^n)$. Given the form
(\ref{paisley}) of each $C^u_{g^n}(\phi), u\in U^{crit}$, that
would then imply that we can determine the sublinear combination
$\Sigma_{u\in U^{crit}} a_u C^u_{g^n}(\phi)$, and the proof of our
third case for the subcase $C_1=0$ would be complete. In order to
determine $\Sigma_{u\in U^{crit}} a_u C^u(g^n)]$, we do the
following:
\newline

\par We may re-express $supdiv_{+}[I^{\frac{n}{2}-A_1}_{g^n}(\phi)]$
in the form:
\begin{equation}
\label{dianaste} \Sigma_{u\in U^{crit}} a_u
Tail_{+}[C^u_{g^n}(\phi)]+ \Sigma_{g\in G} a_g
Tail_{+}[C^g_{g^n}(\phi)]=0
\end{equation}
modulo complete contractions of length $\ge\frac{n}{2}+1$. Here
each $Tail_{+}[C^g_{g^n}(\phi)]$ consists of complete contractions
in the form (\ref{arv6}) and since the sublinear combination
$\Sigma_{g\in G} a_g C^g_{g^n}(\phi)$ is known, we have that the
sublinear combination
\\ $\Sigma_{g\in G} a_g
Tail_{+}[C^g_{g^n}(\phi)]$ is known. Alternatively, in our new
notation using (\ref{clover}):

\begin{equation}
\label{dianaste2} \Sigma_{u\in U^{crit}} a_u
\Sigma_{h=1}^{A_1}C^{u,h}_{g^n}(\phi)+ \Sigma_{g\in G} a_g
Tail_{+}[C^g_{g^n}(\phi)]=0
\end{equation}
modulo complete contractions of length $\ge\frac{n}{2}+1$.  We
will then determine the sublinear combination $\Sigma_{u\in
U^{crit}} a_u C^u(g^n)$ by a trick:
\newline

\par Initially, we polarize the $\frac{n}{2}-A_1$ functions
$\phi$
in the above equation. We denote by $C^{u,h}_{g^n}(\phi_1,\dots
,\phi_{\Delta_1})$ the complete contraction:

\begin{equation}
\label{arv123} contr(R_{i_1j_1k_1l_1}\otimes\dots
\otimes\nabla_{i_1\dots i_{\Delta_1}}R_{i_hj_hk_hl_h}\otimes\dots
\otimes R_{i_{A_1}j_{A_1}k_{A_1}l_{A_1}}\otimes
(\nabla)^{i_1}\phi_1\otimes\dots \otimes
(\nabla)^{i_{\Delta_1}}\phi_{\Delta_1})
\end{equation}

\par We also denote by $\Sigma_{g\in G} a_g Tail_{+}[C^g_{g^n}(\phi_1,\dots ,
\phi_{\Delta_1})]$ the sublinear combination of complete
contractions that arises from $\Sigma_{g\in G} a_g
Tail_{+}[C^g_{g^n}(\phi)]$ by polarizing the $\Delta_1$ functions
$\phi$. It will be a linear combination of complete contractions
in the form:
\begin{equation}
\label{arv800} contr(\nabla^{\Delta_1}_{t_1\dots t_{\Delta_1}}
R_{i_1j_1k_1l_1}\otimes R_{i_2j_2k_2l_2}\otimes \dots \otimes
 R_{i_{A_1}j_{A_1}k_{A_1}l_{A_1}}
\otimes(\nabla)^{i_1}\phi_1\otimes\dots\otimes
(\nabla)^{i_{\frac{n}{2}-A_1}} \phi_{\Delta_1})
\end{equation}
where each $\nabla\phi_h$ contracts against the same factor
$\nabla^{\Delta_1}_{i_1\dots i_{\Delta_1}} R_{i_1j_1k_1l_1}$.
Again, since $\Sigma_{g\in G} a_g C^g_{g^n}(\phi_1, \dots
,\phi_{\Delta_1})$ arises from $\Sigma_{g\in G} a_g
C^g_{g^n}(\phi)$ by polarization, we have
 that the sublinear combination
$\Sigma_{g\in G} a_g C^g_{g^n}(\phi_1,\dots ,\phi_{\Delta_1})$ is
known. Therefore, from (\ref{dianaste2}) we derive an equation
modulo complete contractions of length $\ge\frac{n}{2}+1$:
\begin{equation}
\label{dianaste'} \Sigma_{u\in U^{crit}} a_u
\Sigma_{h=1}^{A_1}C^{u,h}_{g^n}(\phi_1,\dots , \phi_{\Delta_1})]+
\Sigma_{g\in G} a_g Tail_{+}[C^g_{g^n}(\phi_1,\dots ,
\phi_{\Delta_1})]=0
\end{equation}

\begin{definition}
\label{arafat}
 For each $0\le\kappa\le \Delta_1$, we define
$C^{u,h}_{g^n}(\phi_{\kappa +1},\dots ,\phi_{\Delta_1})$ to stand
for the complete contraction:

\begin{equation}
\label{arv123cut}
\begin{split}
& contr(R_{i_1j_1k_1l_1}\otimes\dots
\otimes\nabla_{i_{\kappa +1}\dots
i_{\Delta_1}}R_{i_hj_hk_hl_h}\otimes\dots \otimes
R_{i_{A_1}j_{A_1}k_{A_1}l_{A_1}}\otimes (\nabla)^{i_{\kappa
+1}}\phi_{\kappa +1}
\\& \otimes\dots \otimes (\nabla)^{i_{\Delta_1}}\phi_{\Delta_1})
\end{split}
\end{equation}

It arises from $C^{u,h}_{g^n}(\phi_1,\dots ,
\phi_{\Delta_1})$ by
erasing the factors $\nabla\phi_h, h\le \kappa$ and also erasing
the indices that they contract against in the factor
$\nabla^{\Delta_1}R_{i_hj_hk_hl_h}$. We observe that for $\kappa
=0$, our notation is consistent. We also have for $\kappa
+1=\Delta_1$, we obtain $C^u(g^n)$. We note that by construction
$C^{u,h}_{g^n}(\phi_1,\dots , \phi_{\Delta_1})$ has length
$\frac{n}{2}-\kappa$.
\end{definition}

\par We now consider complete contractions of the form:

\begin{equation}
\label{stein2}
\begin{split}
& contr(\nabla^{\Delta_1-\kappa}_{r_1\dots r_{\Delta_1-\kappa}}
R_{i_1j_1k_1l_1}\otimes R_{i_2j_2k_2l_2}\otimes \dots \otimes
 R_{i_{A_1}j_{A_1}k_{A_1}l_{A_1}}
\otimes(\nabla)^{i_{k+1}}\phi_{\kappa +1}\otimes\dots
\\& \otimes
(\nabla)^{i_{\frac{n}{2}-A_1}} \phi_{\Delta_1})
\end{split}
\end{equation}
where each of the factors $\nabla\phi_h$ contracts against an index in
 the  factor $\nabla^{\Delta_1-\kappa}R_{ijkl}$. We observe that
 up to switching the position of the  factor
 $\nabla^{\Delta_1-\kappa -1}R_{ijkl}$ and a factor
 $R_{i'j'k'l'}$, the complete contractions $C^{u,h}_{g^n}
 (\phi_{\kappa +1},\dots ,\phi_{\Delta_1})$ are in the form
 (\ref{stein2}) above.

\par We now let
$\Sigma_{g\in G^\kappa} a_g C^g_{g^n}(\phi_{\kappa +1},\dots
,\phi_{\Delta_1})$ stand for a generic {\it known} linear
combination of complete contractions in the form (\ref{stein2}).

\par Our claim is then the following:

\begin{lemma}
\label{parakalw} We claim that for any $\kappa,
0\le\kappa\le\Delta_1$, we will have
 that modulo complete contractions of length
$\ge\frac{n}{2}-\kappa +1$:

\begin{equation}
\label{duffy} \Sigma_{u\in U^{crit}} a_u
\Sigma_{h=1}^{A_1}C^{u,h}_{g^n}(\phi_{\kappa +1},\dots
,\phi_{\Delta_1})+ \Sigma_{g\in G^\kappa } a_g
C^g_{g^n}(\phi_{\kappa +1},\dots ,\phi_{\Delta_1})=0
\end{equation}
\end{lemma}

\par Clearly if we can show the above, then using the case
$\kappa =\Delta_1$,
 we will
 then have shown our third case above in the first subcase.
The equation holds exactly because terms
 of greater length have the wrong weight.
\newline

{\it Proof:} We will prove the above by an induction. We assume that we
 know our Lemma for $\kappa=k$ and we will show it for
$\kappa=k +1$, where $k\le \Delta_1$.

We write out our inductive hypothesis:

\begin{equation}
\label{duffy'}
\begin{split}
&L_{g^n}(\phi_{k+1},\dots ,\phi_{\Delta_1})= \Sigma_{u\in
U^{crit}} a_u \Sigma_{h=1}^{A_1}C^{u,h}_{g^n}(\phi_{k+1},\dots
,\phi_{\Delta_1})
\\& +\Sigma_{g\in G^k} a_g C^g_{g^n}(\phi_{k +1},\dots ,
\phi_{\Delta_1})= \Sigma_{y\in Y} a_y C^y_{g^n}(\phi_{k +1},\dots
,\phi_{\Delta_1})
\end{split}
\end{equation}
where each $C^y_{g^n}(\phi_{k +1},\dots ,\phi_{\Delta_1})$ has
length $\ge\frac{n}{2}-k +1$.

\par For each complete contraction
$C_{g^n}(\phi_{k +1},\dots ,\phi_{\Delta_1})$ of weight $-n+2k$ we
define, for the purposes of this proof:

\begin{equation}
\label{stein1} Image^1_{\phi'}[C_{g^n}(\phi_{k +1},\dots
,\phi_{\Delta_1})]=
\partial_{\lambda}|_{\lambda =0}[e^{\lambda(n-2k)\phi'}
C_{e^{2\lambda\phi'}g^n}(\phi_{k +1},\dots ,\phi_{\Delta_1})]
\end{equation}

\par Now, by our inductive hypothesis, we deduce that:

\begin{equation}
\label{rodala}
\begin{split}
&Image^1_{\phi'}\{ \Sigma_{u\in U^{crit}} a_u
\Sigma_{h=1}^{A_1}C^{u,h}_{g^n}(\phi_{k+1},\dots
,\phi_{\Delta_1})\} +
\\&  Image^1_{\phi'}\{ \Sigma_{g\in G^k} a_g
C^g_{g^n}(\phi_{k +1},\dots ,\phi_{\Delta_1})\}=
 Image^1_{\phi'}\{ \Sigma_{y\in Y} a_y
C^y_{g^n}(\phi_{k +1},\dots ,\phi_{\Delta_1})\}
\end{split}
\end{equation}

\par We make a note on how the operation $Image^1_{\phi'}$ acts:
 Consider any complete contraction
$C_{g^n}(\phi_{k+1},\dots ,\phi_{\Delta_1})$ of weight $-n+2k$.
Then, $Image^1_{\phi'}[C_{g^n}(\phi_{k+1},\dots
,\phi_{\Delta_1})]$ is determined as follows: We arbitrarily pick
out one factor $T_{g^n}$ in $C_{g^n}(\phi_{k+1},\dots
,\phi_{\Delta_1})$ and we make all its indices free. We thus have
a tensor $T^{g^n}_{i_1\dots i_h}$. Then, consider all the terms in
$T^{e^{2\phi'}g^n}_{i_1\dots i_h}$ that are linear in $\phi'$ and
involve at least one derivative of $\phi'$. We arbitrarily replace
$T_{g^n}$ in $C_{g^n}(\phi_{k +1},\dots ,\phi_{\Delta_1})$ by one
of those terms, we leave all the other factors unaltered, and
perform the same particular contractions as for
$C_{g^n}(\phi_{k+1},\dots ,\phi_{\Delta_1})$. Adding over all
these arbitrary substitutions, we obtain
$Image^1_{\phi'}[C_{g^n}(\phi_{k+1},\dots ,\phi_{\Delta_1})]$.

\par Now, we restrict our attention to complete contractions
$C_{g^n}(\phi_{k +1},\dots , \phi_{\Delta_1})$ in the form
(\ref{stein2}) and we wish to understand which complete
contractions in
 $Image^1_{\phi'}[C_{g^n}(\phi_{k +1},\dots ,
\phi_{\Delta_1})]$ are in the form:

\begin{equation}
\label{krinos}
\begin{split}
&contr(\nabla^{\Delta_1-\kappa-1}R_{i_1j_1k_1l_1}\otimes
R_{i_2j_2k_2l_2} \otimes
R_{i_{A_1}j_{A_1}k_{A_1}l_{A_1}}\otimes\nabla\phi_{k +2}
\otimes\nabla\phi_{\Delta_1}\otimes \dots\otimes
(\nabla)^h\phi'
\\&\otimes\nabla_h\phi_{k +1})
\end{split}
\end{equation}
In the above complete contraction, the length is $\frac{n}{2}-k
+1$ and each of the factors $\nabla\phi_h, h\ge k+2$ contracts
 against the factor $\nabla^{\Delta_1-k -1}R_{ijkl}$ and the two
 factors $\nabla\phi_{k+1}, \nabla\phi'$ contract between
 themselves. We will call such contractions {\it targets}.
We denote their sublinear combination in each
$Image^1_{\phi'}
[C_{g^n}(\phi_{k+1},\dots , \phi_{\Delta_1})]$ by
$Image^{1,targ}_{\phi'}[C_{g^n}(\phi_{k+1},\dots ,
\phi_{\Delta_1})]$.

\par Now, let us further analyze each
$Image^1_{\phi'}[C_{g^n}(\phi_{k+1},\dots , \phi_{\Delta_1})]$,
where
\\ $C_{g^n}(\phi_{k+1},\dots , \phi_{\Delta_1})$ is in the form
(\ref{stein2}). For each factor $T_f=\nabla^mR_{ijkl}$ ($m\ge 0$,
$1\le f\le A_1$), we denote by
$Full_{T_f}[C_{g^n}(\phi_{k+1},\dots , \phi_{\Delta_1})]$ the sum
of four complete contractions that arises from
$C_{g^n}(\phi_{k+1},\dots , \phi_{\Delta_1})$ by replacing the
factor $T=\nabla^mR_{ijkl}$ by one of the linear expressions
$\nabla^m(\nabla^2\phi'\otimes g)$ on the right hand side of
(\ref{curvtrans}) and then adding those four substitutions. It
follows that each $Full_{T_f}[C_{g^n}(\phi_{k+1},\dots ,
\phi_{\Delta_1})]$ is a sum of four complete contractions of
length $\frac{n}{2}-k$, each in the form:

\begin{equation}
\label{jessica}
\begin{split}
&contr(\nabla^{m_1}R_{ijkl}\otimes\dots
\otimes\nabla^{m_{A_1}-1}R_{i'j'k'l'}\otimes
\nabla^r\phi'\otimes\nabla\phi_{k+1}\otimes\dots\otimes\nabla
\phi_{\Delta_1})
\end{split}
\end{equation}
where $r\ge 2$, and each $m_u\ge 0$. This follows from the
transformation law (\ref{curvtrans}).
\newline

\par On the other hand, for each $C_{g^n}(\phi_{k+1},\dots ,
\phi_{\Delta_1})$, we make note of the one factor
$\nabla^mR_{ijkl}$ with $m>0$ and we call it {\it critical}. We let $LC^{crit}[C_{g^n}(\phi_{k+1},\dots ,
\phi_{\Delta_1})]$ stand for
the sublinear combination that arises in
$Image^1_{\phi'}[C_{g^n}(\phi_{k+1},\dots , \phi_{\Delta_1})]$
when we replace the critical factor by an expression
$\nabla^hR_{ijkl}\nabla^b\phi'$ or $\nabla^hR_{ijkl}\nabla^b\phi'
g_{ab}$, that arises either by virtue of the transformation law
(\ref{levicivita}) or by virtue of the homogeneity of
$R_{ijkl}$ (see (\ref{curvtrans})).

\par Then, for each complete contraction
$C_{g^n}(\phi_{k+1},\dots , \phi_{\Delta_1})$ on the left hand
side of (\ref{rodala}) we have:

\begin{equation}
\label{yanukovich}
\begin{split}
& Image^1_{\phi'}[C_{g^n}(\phi_{k+1},\dots ,
\phi_{\Delta_1})]=\Sigma_{f=1}^{A_1}Full_{T_f}
[C_{g^n}(\phi_{k+1},\dots , \phi_{\Delta_1})]
\\& +LC^{crit}[C_{g^n}(\phi_{k+1},\dots , \phi_{\Delta_1})]
\end{split}
\end{equation}

\par We will now show that:

\begin{equation}
\label{sashka}
\begin{split}& \Sigma_{u\in U^{crit}} a_u
\Sigma_{f=1}^{A_1} Full_{T_f} [C^{u,h}_{g^n}(\phi_{k+1},\dots
,\phi_{\Delta_1})] +
\\& \Sigma_{g\in G} a_g \Sigma_{f=1}^{A_1}
Full_{T_f} [C^g_{g^n}(\phi_{k +2},\dots ,\phi_{\Delta_1})]=
 \Sigma_{j\in J} a_j C^j_{g^n}(\phi_{k+1},\dots ,
\phi_{\Delta_1},\phi')
\end{split}
\end{equation}
where each $\Sigma_{j\in J} a_j C^j_{g^n}(\phi_{k+1},\dots ,
\phi_{\Delta_1},\phi')$ has length $\ge\frac{n}{2}-k +1$ and
is not a target.

\par We see this as follows: Initially, we recall equation
(\ref{rodala}), where the left hand side can be explicitly written
out by virtue of (\ref{yanukovich}) and the right hand side
consists of complete contractions of length
$\ge\frac{n}{2}-k+1$.
This follows from (\ref{curvtrans}) and (\ref{levicivita}).
Therefore, recalling that each $LC^{crit}[C_{g^n}(\phi_{k+1},\dots
,\phi_{\Delta_1})$ in (\ref{yanukovich}) consists of complete
contractions of length $\frac{n}{2}-k+1$, we have:

\begin{equation}
\label{sashka2}
\begin{split}&\Sigma_{u\in U^{crit}} a_u \Sigma_{h=1}^{A_1}\Sigma_{f=1}^{A_1}
Full_{T_f}[C^{u,h}_{g^n}(\phi_{k+1},\dots ,\phi_{\Delta_1})]
\\& +\Sigma_{g\in G^k} a_g \Sigma_{f=1}^{A_1}
Full_{T_f}[C^g_{g^n}(\phi_{k +2},\dots ,\phi_{\Delta_1})]=0
\end{split}
\end{equation}
modulo complete contractions of length $\ge\frac{n}{2}-k +1$.

\par Now, the above holds formally. Hence, there is a sequence of
permutations among the indices of the
 factors in the left hand side of the above with which we can make
the left hand side of the above formally zero, modulo introducing
complete contractions of length $\ge\frac{n}{2}-k +1$. We
want to keep track of the correction terms that arise. We see that
the correction terms can only arise by applying the identity
$[\nabla_A\nabla_B-\nabla_B\nabla_A]X_C=R_{ABCD}X^D$. But we see
that if we apply this identity to a factor $\nabla^mR_{ijkl}$, we
introduce a correction term of length $\frac{n}{2}-k +1$ {\it
which will have a factor $\nabla^r\phi'$, $r\ge 2$}. This is true
because each expression consists of complete contractions in the
form (\ref{jessica}), so there is such a factor to begin with.
Hence, we do not obtain a target in this way.
 On the other hand, if we apply the identity
$[\nabla_A\nabla_B-\nabla_B\nabla_A]X_C=R_{ABCD}X^D$  to the
 factor $\nabla^r\phi'$, $r\ge 2$, we will obtain a correction term
 which will either have a factor $\nabla^u\phi'$, $u\ge 2$ {\it or}  a
factor $\nabla\phi'$ {\it which contracts against a factor
$\nabla^tR_{ijkl}$}. Therefore, we do not obtain a targets in this
way either. We have shown (\ref{sashka}).
\newline

\par Our next claim is:
\newline
\par Claim A: For each $u\in U^{crit}$, $1\le h\le A_1$:

\begin{equation}
\label{byzia2}
\begin{split}
&LC^{crit}[C^{u,h}_{g^n}(\phi_{k+1},\dots , \phi_{\Delta_1})]=
\\&(-2-(\Delta_1 -k-1))
C^{u,h}_{g^n}(\phi_{k+2},\dots ,\phi_{\Delta_1})
(\nabla)^h\phi'\nabla_h\phi_{k+1}
\\& +\Sigma_{j\in J} a_j C^j_{g^n}(\phi_{k+1},\dots
,\phi_{\Delta_1})
\end{split}
\end{equation}
where the linear combination $\Sigma_{j\in J} a_j
C^j_{g^n}(\phi_{k +1},\dots ,\phi_{\Delta_1},\phi')$ is a generic
linear combination of complete contractions of length
$\frac{n}{2}-k+1$ that are not targets.
\newline

\par We show claim A as follows: For each complete contraction $C_{g^n}(\phi_{k+1},\dots
,\phi_{\Delta_1})$ appearing on the left hand side of
(\ref{duffy'}), we have defined
$LC^{crit}[C_{g^n}(\phi_{k+1},\dots ,\phi_{\Delta_1})]$.  Now, we
pay special attention to the one index $i_{k+1}$ in the critical
factor that is contracting against the factor $\nabla\phi_{k+1}$.
Let $LC^{crit,\alpha}[C_{g^n}(\phi_{k+1},\dots
,\phi_{\Delta_1})]$ be the sublinear combination that
 arises in $LC^{crit}[C_{g^n}(\phi_{k+1},\dots
,\phi_{\Delta_1})]$ when we replace the critical factor
$\nabla^{m-k}_{r_{k+1}\dots r_m}R_{ijkl}$ by an expression
$\nabla_{r_{k+1}}\phi'\nabla^{m-k-1}R_{ijkl}$. (Note that the
index $r_{k+1}$ is the one that contracted against the factor
$\nabla\phi_{k+1}$ in $C_{g^n}(\phi_{k+1},\dots
,\phi_{\Delta_1}$)). We denote by
$LC^{crit,\beta}[C_{g^n}(\phi_{k+1},\dots ,\phi_{\Delta_1})]$ the
sublinear combination that arises in
$LC^{crit}[C_{g^n}(\phi_{k+1},\dots ,\phi_{\Delta_1})]$ when we
replace the critical factor in any other way.

\par Hence, $LC^{crit,\beta}[C_{g^n}(\phi_{k+1},\dots
,\phi_{\Delta_1})]$ arises by replacing the critical factor by an
expression in either the form $\nabla^h\phi'\nabla^uR_{ijkl}$,
 $\nabla^h\phi'\nabla^uR_{ijkl}g_{ab}$ with
$h\ge 2$ or of the form $\nabla_\alpha\phi'\nabla^uR_{ijkl}$,
$\nabla_\alpha\phi'\nabla^uR_{ijkl}g_{ab}$
 where the index $\alpha$ is
not the index $r_{k+1}$ that contracts against $\nabla\phi_{k
+1}$.

\par We  observe that the sublinear combinations
$LC^{crit,\alpha}[C^{u,h}_{g^n}(\phi_{k+1},\dots
,\phi_{\Delta_1})]$, $LC^{crit,\alpha}[C^g_{g^n}(\phi_{k+1},\dots
,\phi_{\Delta_1})]$ consist of targets, whereas the sublinear
combinations $LC^{crit,\beta}[C^{u,h}_{g^n}(\phi_{k+1},\dots
,\phi_{\Delta_1})]$, $LC^{crit,\beta}[C^g_{g^n}(\phi_{1+1},\dots
,\phi_{\Delta_1})]$ contain no targets.

\par  Therefore, in view of the above, in order to show Claim A,
we only have to show that for each $u\in U^{crit}$ and each $1\le
h\le A_1$, we have that:

\begin{equation}
\label{elaine}
\begin{split}
&LC^{crit,\alpha}[C^{u,h}_{g^n}(\phi_{k+1},\dots
,\phi_{\Delta_1})]=
\\& (-2-(\Delta_1-k -1))\cdot C^{u,h}_{g^n}
(\phi_{k+2},\dots ,\phi_{\Delta_1})(\nabla)^h\phi_1\nabla_h\phi'
\end{split}
\end{equation}

\par Hence, we only have to show that the sublinear combination of
expressions in $Image^1_{\phi'}
[\nabla^{\Delta_1-k}_{r_{k +1}\dots r_{\Delta_1}}R_{ijkl}]$ that
are in the form $\nabla_{r_{k+1}}\phi'
\nabla^{\Delta_1-k}_{r_{k+2}\dots r_{\Delta_1-\kappa}}
R_{ijkl}$ is
precisely $(-2-(\Delta_1-k -1))\cdot \nabla_{r_{k+1}}\phi'
\nabla^{\Delta_1-\kappa}_{r_{k+2}\dots
r_{\Delta_1-\kappa}}R_{ijkl}$.
 But this is only a matter of applying (\ref{levicivita}) to all
  the pairs $(r_{k+1}, r_a)$, $a\ge k+2$ and the pairs
$(r_{k+1},i),\dots (r_{k+1},l)$
  and also by taking into account the expression
\\ $2\nabla_{r_{k+1}}\phi'\nabla^{\Delta_1-k -1}_{r_{k+2}\dots
r_{\Delta_1}} R_{ijkl}$ that arises by virtue of the homogeneity
of the factor $R_{ijkl}$.

Combining the equations (\ref{rodala}), (\ref{yanukovich}),
(\ref{sashka}), (\ref{byzia2}) and (\ref{elaine}) above, we have
that:

\begin{equation}
\label{neeo}
\begin{split}
&\Sigma_{j\in J} a_j C^j_{g^n}(\phi_{k+1},\dots
,\psi_{\Delta_1},\phi') +\Sigma_{u\in U^{crit}}
a_u\Sigma_{h=1}^{A_1}
LC^{crit,\alpha}[C^{u,h}_{g^n}(\phi_{k+1},\dots ,\phi_{\Delta_1})]
\\& +\Sigma_{g\in G^k} a_g
LC^{crit,\alpha}[C^g_{g^n}(\phi_{k+1},\dots ,\phi_{\Delta_1})]
 +\Sigma_{u\in U^{crit}}
a_u\Sigma_{h=1}^{A_1}
LC^{crit,\beta}[C^{u,h}_{g^n}(\phi_{k+1},\dots ,\phi_{\Delta_1})]
\\& +\Sigma_{g\in G^k} a_g
LC^{crit,\beta}[C^g_{g^n}(\phi_{k+1},\dots ,\phi_{\Delta_1})] =
\Sigma_{z\in Z} a_z C^z_{g^n}(\phi_1,\dots ,\phi_{\Delta_1},\phi')
\end{split}
\end{equation}
where the sublinear combination $\Sigma_{z\in Z} a_z
C^z_{g^n}(\phi_{k+1},\dots ,\phi_{\Delta_1},\phi')$ stands for:

$$Image^1_{\phi'}[\Sigma_{y\in Y} a_y C^y_{g^n}(\phi_{k+1},\dots
,\phi_{\Delta_1})],$$ and hence each $C^z_{g^n}(\phi_{k+1},\dots
,\phi_{\Delta_1},\phi')$ either has length $\ge \frac{n}{2}-\kappa
+2$ or has length $\ge \frac{n}{2}-k +1$ but has a factor
$\nabla^u\phi', u\ge 2$ (so it is not a target). Therefore, since
(\ref{neeo}) must hold formally, we deduce that, modulo complete
contractions of length $\ge\frac{n}{2}-k +2$:

\begin{equation}
\label{shila}
\begin{split}
&\Sigma_{u\in U^{crit}} a_u \Sigma_{h=1}^{A_1}
LC^{crit,\alpha}[C^{u,h}_{g^n}(\phi_{k+1},\dots ,\phi_{\Delta_1})]
\\& +\Sigma_{g\in G^k} a_g
LC^{crit,\alpha}[C^g_{g^n}(\phi_{k+1},\dots ,\phi_{\Delta_1})]=0
\end{split}
\end{equation}

\par Now, since we are assuming that the linear combination
$\Sigma_{g\in G^k} a_g C^g_{g^n}(\phi_{k+1},\dots
,\phi_{\Delta_1})$ is known, we deduce that the linear combination
$\Sigma_{g\in G^k} a_g LC^{crit,\alpha}[C^g_{g^n}(\phi_{k+1},\dots
,\phi_{\Delta_1})]$ is also known.

\par  We have thus completed the proof of the third case if $C_1=0$.
$\Box$
\newline

{\it The subcase $C_1>0$:}
\newline

\par The second subcase is almost entirely similar. We again write
out $I^{\frac{n}{2}}_{g^n}(\phi)$ in the form
(\ref{ishbilia}). We
write $C_1=\gamma$. We consider the Shadow divergence formula for
$I^{\frac{n}{2}-A_1}_{g^n}(\phi)$ and we focus on the sublinear
 combination $Shad_{+}[I^{\frac{n}{2}-A_1}_{g^n}(\phi)]$
in $Shad[I^{\frac{n}{2}-A_1}_{g^n}(\phi)]$ which consists of
$\vec{\xi}$-contractions in the form:

\begin{equation}
\label{arv8} \begin{split} &contr(\nabla^{\Delta_1}_{t_1\dots
t_{\Delta_1}} R_{i_1j_1k_1l_1}\otimes\dots \otimes
R_{i_2j_2k_2l_2}\otimes
R_{i_{A_1-\gamma}j_{A_1-\gamma}k_{A_1-\gamma}l_{A_1-\gamma}})\\& \otimes
(|\vec{\xi}|^2)^\gamma\otimes
(\nabla)^{s_1}\phi\otimes\dots\otimes
(\nabla)^{s_{\frac{n}{2}-A_1}}\phi)
\end{split}
\end{equation}
where each a factor $\nabla\phi$ contracts against an index in the
factor $\nabla^{\Delta_1}R_{ijkl}$.

\par As in all the previous cases, we have that:

\begin{equation}
\label{mcmillan}
Shad_{+}[I^{\frac{n}{2}-A_1}_{g^n}(\phi)]=0
\end{equation}
modulo $\vec{\xi}$-contractions of length $\ge\frac{n}{2}+1$,
since the shadow divergence formula holds formally.

\par As before, for each $C^u_{g^n}(\phi)$ that is subsequent to the
 critical list, we have that $Tail^{Shad}_{+}[C^u_{g^n}(\phi)]=0$.
Hence, we have that:

\begin{equation}
\label{mcmillan2} \Sigma_{g\in G} a_g Shad_{+}[C^g_{g^n}(\phi)]+
\Sigma_{u\in U^{crit}} a_u Shad_{+}[C^u_{g^n}(\phi)]=0
\end{equation}
modulo $\vec{\xi}$-contractions of $\vec{\xi}$-length
$\ge\frac{n}{2}+1$.

\par Moreover, for each $u\in U^{crit}$, where
$C^u_{g^n}(\phi)$ is in the form:

\begin{equation}
\label{arv7} contr(R_{i_1j_1k_1l_1}\otimes\dots \otimes
R_{i_{A_1-\gamma}j_{A_1-\gamma}k_{A_1-\gamma}l_{A_1-\gamma}}\otimes
R^\gamma\otimes \Delta\phi\otimes\dots\otimes\Delta\phi)
\end{equation}

we have that $Tail^{Shad}_{+}[C^u_{g^n}(\phi)]$ can be written out
as:

\begin{equation}
\label{persia} Tail^{Shad}_{+}[C^u_{g^n}(\phi)]=
\Sigma_{h=1}^{A_1-\gamma}C^{u,h}_{g^n}(\phi)
\end{equation}

where $C^{u,h}_{g^n}(\phi)$ is in the form:

\begin{equation}
\label{hobsbaum} \begin{split} &contr(R_{i_1j_1k_1l_1}\otimes\dots
\otimes \nabla^{i_1\dots
i_{\frac{n}{2}-A_1}}R_{i_hj_hk_hl_h}\otimes\dots\otimes
R_{i_{A_1-\gamma}j_{A_1-\gamma}k_{A_1-\gamma}l_{A_1-\gamma}}
\\& \otimes
(|\vec{\xi}|^2)^\gamma
\otimes\nabla_{i_1}\phi\otimes\dots\otimes
\nabla_{i_{\frac{n}{2}-A_1}}\phi)
\end{split}
\end{equation}

\par We then define $C^u(g^n)$ to stand for the complete contraction:

$$contr(R_{i_1j_1k_1l_1}\otimes\dots \otimes
R_{i_{A_1-\gamma}j_{A_1-\gamma}k_{A_1-\gamma}l_{A_1-\gamma}})$$

\par Hence, using the equation (\ref{mcmillan2}) and repeating the same
 argument as in the above case, we may determine the
 sublinear combination $\Sigma_{u\in U^{crit}} a_u
C^u(g^n)$, and hence also the sublinear combination $\Sigma_{u\in
U^{crit}} a_u C^u_{g^n}(\phi)$. We have completed the proof of
Lemma \ref{A>0}. $\Box$
\newline

\end{document}